\documentclass[11pt,leqno,english]{article}

\usepackage{ShVa718_719}

\newcommand{\EquationRefCtrl}{Part}

\newcommand{\MaybeMinus}[2]{%
	\IfEqual{#1}{-}%
		{}%
		{#2}%
}

\newcommand{\IfMinus}[3]{%
	\IfEqual{#1}{-}%
		{#2}%
		{#3}%
}

\newcommand{\Her}[1]{\MathX{%
	H(#1)%
}}

\newcommand{\HerStr}[3][]{\MathX{%
	\Structure[#1]{\Her{#2}, \in, #3}%
}}

\newcommand{\One}{\MathX{%
	\mathbf{1}%
}}

\newcommand{\Name}[2][]{\MathX{%
	\IfEmpty{#1}%
		{\tilde{#2}}%
		{\widetilde{#2}}%
}}

\newcommand{\Subsets}[2]{\MathX{%
	{[#1]}^{#2}%
}}

\renewcommand{\BigUnion}[2]{\MathX{%
	\IfMinus{#1}%
		{\bigcup #2}%
		{\bigcup_{#1} #2}%
}}

\newcommand{\Product}[2]{\MathX{%
	\IfMinus{#1}%
		{\prod #2}%
		{\prod_{#1} #2}%
}}

\newcommand{\br}[4][b]{\MathX{%
	\IfMinus{#3}%
		{{#1}^{#2}_{#4}}%
		{{#1}^{#2}_{#3,#4}}%
}}

\newcommand{\Namebr}[3]{%
	\br[\Name{b}]{#1}{#2}{#3}%
}

\newcommand{\Ind}[2]{\MathX{%
	\IfMinus{#2}%
		{{\Delta}^{#1}}%
		{{\Delta}^{#1}_{#2}}%
}}

\newcommand{\Tree}[2]{\MathX{%
	\IfMinus{#2}%
		{{T}^{#1}}%
		{{T}^{#1}_{#2}}%
}}

\newcommand{\Br}[2]{\MathX{%
	\Symb{Br}_{#1}(#2)%
}}

\newcommand{\MuSeq}{\MathX{%
	\bar{\mu}%
}}

\newcommand{\ChiSeq}[1][]{%
	\MathX{\bar{\chi}%
}}

\newcommand{\YSeq}[1][]{\MathX{%
	\bar{y}%
}}

\newcommand{\ZSeq}[1][]{\MathX{%
	\bar{z}%
}}

\newcommand{\ClosedCardSet}{\MathX{%
	{\Omega}_{\MuSeq}%
}}

\newcommand{\Fn} [2] {\MathX{%
	\Symb{Fn}(#1,2,#2)%
}}

\newcommand{\EquivSymbol}{\thicksim}

\newcommand{\Equiv}[1][\phi,\Param]{\MathX{%
	\mathrel{{\EquivSymbol}_{\scriptscriptstyle #1}}%
}}

\newcommand{\NotEquiv}[1][\phi,\Param]{%
	\mathrel{{\not\EquivSymbol}_{\scriptscriptstyle #1}}%
}

\newcommand{\EquivName}[1][\phi,\ParamName]{\MathX{%
	\mathrel{{\EquivSymbol}_{\scriptscriptstyle #1}}%
}}

\newcommand{\NotEquivName}[1][\phi,\ParamName]{%
	\mathrel{{\not\EquivSymbol}_{\scriptscriptstyle #1}}%
}

\newcommand{\GenExt}[1]{\MathX{%
	V[#1]%
}}

\newcommand{\OfGenExt}[3][]{\MathX{%
	\Par[#1]{#2}^{\GenExt{#3}}%
}}

\newcommand{\GenExtSat}[2]{\MathX{%
	V[#1] \models #2%
}}

\newcommand{\DoubleGenExt}[2]{\MathX{%
	V[#1][#2]%
}}

\newcommand{\DoubleGenExtSat}[3]{\MathX{%
	V[#1][#2] \models #3%
}}

\newcommand{\InExt}[3][]{\MathX{%
	#1(#2 {#1)}^{V[#3]}%
}}

\newcommand{\Param}{\MathX{R}}

\newcommand{\ParamName}{\MathX{%
	\Name{R}%
}}

\newcommand{\No}[1]{\MathX{
	\Symb{No}(#1)%
}}

\newcommand{\Forcing}[1]{\MathX{%
	\mathbb{P}(#1)%
}}

\newcommand{\Forces}[2][]{%
	\IfMinus{#1}%
		{\Vdash_{#2}}%
		{\Vdash_{\Forcing {#2}}}%
}

\newcommand{\FunName}[1]{\MathX{%
	{\Name{\mathcal{F}}}_{#1}%
}}

\newcommand{\Maps}[1]{\MathX{%
	\Symb{Mps}(#1)%
}}

\newcommand{\Map}[4][]{\MathX{%
	\IfEmpty{#1}{%
		\def\pvtempa{#2}%
	}{
		\def\pvtempa{{#2}^{#1}}%
	}%
	\IfEqual{#3}{1}{%
		\def\pvtempb{\pvtempa_{1st}}%
	}{
		\IfMinus{#3}{%
			\def\pvtempb{\pvtempa}%
		}{
			\def\pvtempb{\pvtempa_{#3}}%
		}%
	}%
	\IfMinus{#4}{%
		\pvtempb%
	}{
		\pvtempb(#4)%
	}%
}}

\newcommand{\Autom}[1]{\MathX{%
	\hat{#1}%
}}

\newcommand{\Ps}[2]{\MathX{%
	\IfMinus{#2}{%
		\Delta^{#1}%
	}{
		\IfEqual{#2}{1}%
			{\Delta_{1st}^{#1}}%
			{\Delta_{#2}^{#1}}%
	}%
}}

\newcommand{\SEp}[2]{\MathX{%
	\Pair {\xi^{#1}_{#2}} {\delta^{#1}_{#2}}%
}}

\newcommand{\IClass}[1][*]{\MathX{%
	\Lambda^{#1}%
}}

\newcommand{\XiSeq}[1]{\MathX{%
	{\bar{\xi}}^{#1}%
}}

\newcommand{\DeltaSeq}[1]{\MathX{%
	\IfMinus{#1}{%
		\bar{\delta}%
	}{%
		{\bar{\delta}}^{#1}%
	}%
}}

\newcommand{\EpsilonSeq}[1]{\MathX{%
	\IfMinus{#1}{%
		\bar{\epsilon}%
	}{%
		{\bar{\epsilon}}^{#1}%
	}%
}}

\newcommand{\MuCard}[1]{\MathX{%
	\mu_{\xi^*_{#1}}%
}}

\newcommand{\EName}[1]{\MathX{%
	\tau_{#1}%
}}

\newcommand{\EquivENames}[2]{\MathX{%
	\EName{#1} \EquivName \EName{#2}%
}}

\newcommand{\NotEquivENames}[2]{\MathX{%
	\EName{#1} \NotEquivName \EName{#2}%
}}

\newcommand{\ESeqs}{\MathX{%
	{\mathcal{E}}^{*}%
}}

\newcommand{\EpsilonPair}[2]{\MathX{%
	\IfMinus{#1}{%
		\Pair{\xi^{*}_{#2}}{\epsilon_{#2}}%
	}{
		\Pair {\xi^{*}_{#2}} {\epsilon^{#1}_{#2}}%
	}%
}}

\newcommand{\DeltaPair}[2]{\MathX{%
	\IfMinus{#1}{%
		\Pair{\xi^{*}_{#2}}{\delta_{#2}}%
	}{
		\Pair {\xi^{*}_{#2}} {\delta^{#1}_{#2}}%
	}%
}}

\newcommand{\Agr}[2]{\MathX{%
	A(#1,#2)%
}}

\newcommand{\CritInd}[1][\alpha^*]{\MathX{%
	I^{#1}%
}}

\newcommand{\CritJ}{\MathX{%
	J^*%
}}

\newcommand{\Repr}{\MathX{%
	\Gamma%
}}

\renewcommand{\epsilon}{\MathX{%
	\varepsilon%
}}

\newcommand{\SigmaAgr}{\MathX{%
	\mathcal{A}^*%
}}

\newcommand{\SmallCards}{\MathX{%
	S(\ParamName)%
}}

\begin{document}

%
\title{\protect %
On equivalence relations $\Sigma_1^1$-definable over \Her \kappa
}

\author{%
	Saharon Shelah%
	\thanks{%
	Research supported by the United States-Israel
	Binational Science Foundation. Publication 719.%
	}\\%
	\and%
	Pauli V\"{a}is\"{a}nen%
}

\date{November 12, 1999}
\maketitle

\begin{abstract}%
Let $\kappa$ be an uncountable regular cardinal.  Call an equivalence
relation on functions from $\kappa$ into 2 $\Sigma_1^1$-definable over
\Her \kappa if there is a first order sentence $\phi$ and a parameter
$\Param \Subset \Her \kappa$ such that functions $f ,g \in \Functions
\kappa 2$ are equivalent iff for some $h \in \Functions \kappa 2$, the
structure \HerStr \kappa {\Param, f, g, h} satisfies $\phi$, where
$\in$, \Param, $f$, $g$, and $h$ are interpretations of the symbols
appearing in $\phi$. All the values $\mu$, $1 \leq \mu \leq \SuccCard
\kappa$ or $\mu = 2^\kappa$, are possible numbers of equivalence
classes for such a $\Sigma_1^1$-equivalence relation. Additionally,
the possibilities are closed under unions of $\leq \kappa$-many
cardinals and products of $< \kappa$-many cardinals. We prove that,
consistent wise, these are the only restrictions under the singular
cardinal hypothesis. The result is that the possible numbers of
equivalence classes of $\Sigma_1^1$-equivalence relations might
consistent wise be exactly those cardinals which are in a prearranged
set, provided that the singular cardinal hypothesis holds and that the
following necessary conditions are fulfilled: the prearranged set
contains all the nonzero cardinals in $\kappa \Union \Braces {\kappa,
\SuccCard \kappa, 2^\kappa}$ and it is closed under unions of $\leq
\kappa$-many cardinals and products of $< \kappa$-many cardinals.  The
result is applied in \cite {ShVaWeakly} to get a complete solution of
the problem of the possible numbers of strongly equivalent
non-isomorphic models of weakly compact cardinality.
\footnote{%
}%
\end{abstract}
%

\begin{SECTION} {-} {Introduction} {Introduction}
%
%
We deal with equivalence relations which are in a simple way definable
over \Her \kappa when $\kappa$ is an uncountable regular cardinal. The
conclusion will be that we can completely control the possible numbers
of equivalence classes of such equivalence relations, provided that
the singular cardinal hypothesis holds. The main application of this
is the solution of the problem of the possible numbers of strongly
equivalent non-isomorphic models of weakly compact
cardinality. Namely, we prove in \cite {ShVaWeakly} that when $\kappa$
is a weakly compact cardinal, there exists a model of cardinality
$\kappa$ with $\mu$-many strongly equivalent non-isomorphic models if,
and only if, there exists an equivalence relation which is
$\Sigma_1^1$-definable over \Her \kappa and it has exactly $\mu$
different equivalence classes.  The paper \cite {ShVaWeakly} can be
read independently of this paper, if the reader accepts the present
conclusion on faith.

For every nonzero cardinals $\mu \leq \kappa$ or $\mu = 2^\kappa$,
there is an equivalence relation $\Sigma_1^1$-definable over \Her
\kappa with $\mu$ different equivalence classes.  There is also a
$\Sigma_1^1$-equivalence relation with \SuccCard \kappa-many classes
\Note {\Lemma {Tree_and_Eq}}. Furthermore, by a simple coding, the
possible numbers of equivalence classes of $\Sigma_1^1$-equivalence
relations are closed under unions of length $\leq \kappa$ and products
of length $< \kappa$. In other words, assuming that $\gamma \leq
\kappa$ and $\chi_i$, $i < \gamma$, are cardinals such that for each
$i < \gamma$, there is a $\Sigma_1^1$-equivalence relation having
$\chi_i$ different equivalence classes, then there is a
$\Sigma_1^1$-equivalence relation having \BigUnion {i < \gamma}
{\chi_i} different equivalence classes, and if $\gamma < \kappa$,
there is also a $\Sigma_1^1$-equivalence relation with \Card {\Product
{i < \gamma} {\chi_i}} different equivalence classes \Note {\Lemma
{Eq_closure}}.

What are the possible numbers of equivalence classes between \SuccCard
\kappa and $2^\kappa$?  The existence of a tree $T \Subset \Her
\kappa$ with $\mu$ different $\kappa$-branches through it implies that
there is a $\Sigma^1_1$-equivalence relation having exactly $\mu$
equivalence classes \Note {\Lemma {Tree_and_Eq}}. Therefore, existence
of a Kurepa tree of height $\kappa$ with more than \SuccCard
\kappa-many and less than $2^\kappa$-many $\kappa$-branches through it
presents an example of a $\Sigma^1_1$-equivalence relation with many
equivalence classes, but not the maximal number.
On the other hand, in an ordinary Cohen extension of $L$, in which
$2^\kappa > \SuccCard \kappa$, there is no definable equivalence
relation having $\mu$-many different equivalence classes when
$\SuccCard \kappa < \mu < 2^\kappa$ \Note {a proof of this fact is
straightforward, and in fact, it is involved in the proof presented in
\Section {Sigma11No}}.

We show that, consistent wise, the closure properties mentioned are
the only restrictions concerning the possible numbers of equivalence
classes of $\Sigma^1_1$-equivalence relations. Namely the conclusion
will be the following:
Suppose $\lambda > \SuccCard \kappa$ is a cardinal with
$\lambda^\kappa = \lambda$ and $\Omega$ is a set of cardinals between
\SuccCard \kappa and $\lambda$ so that it is closed under unions of
$\leq \kappa$-many cardinals and products of $< \kappa$-many
cardinals. We shall prove that after adding into $L$ in the
``standard'' way Kurepa trees of height $\kappa$ with $\mu$-many
$\kappa$-branches through it, for every $\mu \in \Omega$ \Note {and
repeating each addition $\lambda$-many times}, there exists, in the
generic extension, an equivalence relation $\Sigma_1^1$-definable over
\Her \kappa with $\mu$-many equivalence classes if, and only if, $\mu$
is a nonzero cardinal $\leq \SuccCard \kappa$ or $\mu$ is in $\Omega
\Union {\Braces {2^\kappa}}$.

In order to make this paper self contained, we introduce the standard
way to add a Kurepa tree and give some basic facts concerning that
forcing in \Section {KurepaTree}. The essential points are the
following. Firstly, if one adds several new Kurepa trees, the addition
of new trees does not produce new $\kappa$-branches of the old
trees. Secondly, permutations of ``the labels'' of the
$\kappa$-branches of the generic Kurepa trees, determine many
different automorphisms of the forcing itself. These kind of
automorphisms can be used ``to copy'' two different equivalence
classes of a definable equivalence relation to several different
equivalence classes. In fact, this is ``the straightforward way'' to
show that in an ordinary Cohen extension of $L$, a definable
equivalence relation has either $\leq \SuccCard \kappa$-many, or the
maximal number $2^\kappa$-many equivalence classes.
The main difference, however, between the standard ``Cohen-case'' and
the proof presented in \Section {Sigma11No} is that the ordinary
$\Delta$-lemma cannot be directly applied as can be done in the former
case.

In \Section {Sigma11} we introduce proofs of some basic facts
mentioned above. The crucial fact is that a $\Sigma^1_1$-equivalence
relation is absolute for various generic extensions \Note {\Lemma
{Absolute} and \Conclusion {Absolute}}. The theorem is formally
written in \Section {Sigma11No}, and the proof of it is divided into
several subsections. The main idea is the following.
We start to look at an equivalence relation which is
$\Sigma^1_1$-definable over \Her \kappa using some parameter of
cardinality $\kappa$.
The forcing consist of addition of $\lambda$-many different
trees. However, we may assume that the forcing name of the parameter
has cardinality $\kappa$, and thus, there are only $\kappa$-many trees
which really has ``effect'' on the number of classes of the fixed
equivalence relation.  So we restrict ourselves to the subforcing
consisting of the addition of these $\kappa$-many ``critical''
trees. \Note {Note, in \Lemma {Subforcing} we introduce a subforcing
consisting of addition of \SuccCard \kappa-many trees, but right after
that in Subsection \ref {isom_names}, we define ``isomorphism
classes'' of names in order to concentrate only on $\kappa$-many
generic trees.}
Then as explained in Subsection \ref {small}, from our assumption that
the singular cardinal hypothesis holds, it immediately follows that
either 1) the fixed equivalence relation has $\chi$ classes, where
$\chi$ is a union of $\leq \kappa$-many cardinals or a product of $<
\kappa$-many cardinals from the prearranged set $\Omega$, or
otherwise, 2) the number of equivalence classes really depends on
$\kappa$-many trees, not less than $\kappa$-many.
On the other hand, we know that the rest of the forcing, i.e., the
addition of the other trees than those $\kappa$-many critical ones,
produces $\lambda$-many new subsets of any set having size
$\kappa$. So, when the equivalence depends on $\kappa$-many trees, we
show in Subsection \ref {large} that either 1) the fixed equivalence
relation has $\chi$ classes, where $\chi$ is a union of $\kappa$-many
products having length $< \kappa$ and cardinals from $\Omega$, or
otherwise, 2) the rest of the forcing produces $\lambda$-many new
equivalence classes.

In \Section {Remarks} we give some concluding remarks.
%
%
\end{SECTION}

\begin{SECTION} {-} {KurepaTree} {Adding Kurepa trees}
%
%
Throughout of this paper we assume that $\kappa$ is an uncountable
regular cardinal and $\kappa^{<\kappa} = \kappa$.
For sets $X$ and $Y$ we denote the set of all functions from $X$ into
$Y$ by \Functions X Y. For a cardinal $\mu$, we let \Subsets X \mu be
the set of all subsets of $X$ having cardinality $\mu$.

The following forcing is the ``standard'' way to add a Kurepa tree
\cite {Jech71,JechSetTheory}.

\begin{DEFINITION}{TForc}%
\Lapikayty{1.10.}%
Let $\mu$ be a cardinal $\geq \kappa$. Define a forcing $P_\mu$ as
follows. It consists of all pairs $p = \Pair {\Tree p -} {\Seq {\br p
- \delta} {\delta \in \Ps p -}}$ where 
\begin{itemize}

\item for some $\alpha < \kappa$, \Tree p - is a subset of \Set {\eta}
{\eta \in \Functions \beta 2 \And \beta < \alpha} such that it is of
cardinality $< \kappa$ and closed under restriction;

\item \Ps p - is a subset of $\mu$ having cardinality $<\kappa$ and
each \br p - \delta is an $\alpha$-branch trough \Tree p - when \Tree
p - is ordered by the inclusion.

\end{itemize}
For all $p, q \in P_\mu$, we define that $q \leq p$ if 
\begin{itemize}

\item \Tree q - is an end-extension of \Tree p -;

\item $\Ps p - \Subset \Ps q -$;

\item for every $\delta \in \Ps p -$, \br q -\delta is an extension of
\br p -\delta.

\end{itemize}
\end{DEFINITION}

\begin{FACT}{basic}%
\begin{ITEMS}%
\item $P_\mu$ is $\kappa$-closed and it satisfies \SuccCard
\kappa-chain condition.

\ITEM{branches}%
Suppose $G$ is a $P_\mu$-generic set over $V$.  In \GenExt G, $\Tree G
- = \BigUnion {p \in G} {\Tree p -}$ is a tree of height $\kappa$ and
each of its level has cardinality $< \kappa$.

\end{ITEMS}
\end{FACT}

\begin{LEMMA}{no_new_branches}%
\Lapikayty{1.10.}%
Let \Name Q be such that $\One \Forces [-] {P_\mu} \Quote {\Name Q
\Text {is a $\kappa$-closed forcing notion}}$. Suppose $G$ is a
$P_\mu$-generic set over $V$ and $H$ is $Q$-generic set over \GenExt
G.  Then, in \DoubleGenExt G H, the $\kappa$-branches trough the tree
$\Tree G - = \BigUnion {p \in G} {\Tree p -}$ are the functions \br G
-\delta, $\delta < \mu$, having domain $\kappa$ and satisfying for
every $\alpha < \kappa$ that $\br G - \delta (\alpha) = \br p - \delta
(\alpha)$ for some $p \in G$ with $\delta \in \Ps p -$ and $\alpha \in
\Dom {\br p - \delta}$.

\begin{Proof}%
The idea of the proof is the same as in \cite {Jech71}. Suppose \Pair
{p_0} {{\Name q}_0} is a condition in $P_\mu * \Name Q$ and \Name t is
a name such that
\[
	\Pair {p_0} {{\Name q}_0} \Forces [-] {P_\mu * \Name Q}
	\Quote{
	\Name t \Text{is a $\kappa$-branch through}
	\Name {\Tree G -}
	\And
	\Name t \not\in \Set {\Namebr G -\delta}{\delta < \mu}}.
\]
Since $\One \Forces [-] {P_\mu * \Name Q} \Quote{\kappa \Text{is a
regular cardinal}}$, it follows that every condition below $\Pair
{p_0} {{\Name q}_0}$ forces that for all $X \in \Subsets \mu {<
\kappa}$ and $\beta < \kappa$, there is $\alpha > \beta$ with
\[
	\Name t (\alpha) \not\in
	\Set {\Namebr G -\delta (\alpha)}{\delta \in X}.
\]
Let $\alpha_0$ be the height of \Tree {p_0} -. Choose conditions \Pair
{p_n} {{\Name q}_n} from $P_\mu * \Name Q$ and ordinals $\alpha_n$, $1
< n < \omega$, so that for every $n < \omega$, the height of the tree
\Tree {p_{n+1}} - is greater than $\alpha_n$, $\Pair {p_{n+1}} {{\Name
q}_{n+1}} \leq \Pair {p_n} {{\Name q}_n}$, and
\EQUATION{new_branch}{%
	\Pair {p_{n+1}} {{\Name q}_{n+1}}
	\Forces [-] {(P_\mu * \Name Q)}
	\,\Name t (\alpha_{n+1})
	\not\in \Set {\Namebr G -\delta (\alpha_{n+1})}
	{\delta \in \Ps {p_n} -}.
}%
Define $r$ to be the condition in $P_\mu$ satisfying $\Tree r - =
\BigUnion {n < \omega} {\Tree {p_n} -}$, $\Ps r - = \BigUnion {n <
\omega} {\Ps {p_n} -}$, and for every $\delta \in \Ps r -$, $\br r -
\delta = \BigUnion {n \in (\omega \Minus m)} {\br {p_n} - \delta}$,
where $m$ is the smallest index with $\delta \in \Ps {p_m} -$.  Then
\Tree r - is of height $\alpha = \BigUnion {n < \omega} {\alpha_n}$.
In order to restrict the \Th \alpha level of the generic tree,
abbreviate the function \BigUnion {\gamma < \alpha} {\br r - \delta
(\gamma)}, $\delta \in \Ps r -$, by $f_\delta$, and define $r'$ to be
the condition in $P_\mu$ with $\Tree {r'} - = \Tree r - \Union \Set
{f_\delta} {\delta \in \Ps r -}$, $\Ps {r'} - = \Ps r -$, and for
every $\delta \in \Ps {r'} -$ and $\beta \leq \alpha$,
\[
	\FunctionDefinition{\br {r'} - \delta (\beta)}{
		\FunctionDefMidCase
		{\br r - \delta (\beta)} {\If \beta < \alpha}
		\FunctionDefLastCase
		{f_\delta} {\If \beta = \alpha}
	}
\]
Now $r'$ forces that the \Th \alpha level of the generic tree \Name
{\Tree G -} consist of the elements $f_\delta$, $\delta \in \Ps {r'}
-$.

Since $r'$ forces \Name Q to be $\kappa$-closed and \Seq {{\Name q}_n}
{n < \omega} to be a decreasing sequence of conditions, there is
$\Name q'$ so that $\Pair {r'} {q'} \leq \Pair {p_n} {{\Name q}_n}$
for every $n < \omega$. Since \Pair {r'} {q'} forces that $\Name t
(\alpha) \in \Set {f_\delta} {\delta \in \Ps {r'} -}$, there are
$\delta \in \Ps {r'} -$ and a condition $\Pair {r''} {\Name q''} \leq
\Pair {r'} {\Name q'}$ in $P_\mu * \Name Q$ forcing that $\Name t
(\alpha) = f_\delta$. However, if $n$ is the smallest index with
$\delta \in \Ps {p_n} -$, then \Pair {r''} {\Name q''} forces that
\[
	\Name t (\alpha_{n+1})
	= \Res {f_\delta} {\alpha_{n+1}}
	= \br {r'} - \delta (\alpha_{n+1})
	= \br {r''} - \delta (\alpha_{n+1})
	= \Namebr G - \delta (\alpha_{n+1}),
\]
contrary to \Equation {new_branch}.
\end{Proof}%
\end{LEMMA}



\begin{DEFINITION}{ProdForc}%
\Lapikayty{1.10.}%
Suppose $\lambda > \SuccCard \kappa$ is a cardinal with
$\lambda^\kappa = \lambda$. Let $\MuSeq = \Seq {\mu_\xi} {\xi <
\lambda}$ be a fixed sequence of cardinals such that $\kappa < \mu_\xi
\leq \lambda$ and for every $\chi \in \Set {\mu_\xi} {\xi < \lambda}
\Union \Braces {\lambda}$, the set \Set {\zeta < \lambda} {\mu_\zeta =
\chi} has cardinality $\lambda$.
We define \Forcing \MuSeq to be the product of $P_{\mu_\xi}$ forcings:
\begin{itemize}

\item \Forcing \MuSeq is the set of all functions $p$ such that \Dom p
is a subset of $\lambda$ with cardinality $< \kappa$, and for every
$\xi \in \Dom p$, $p(\xi)$ is a condition in $P_{\mu_\xi}$;

\item the order of \Forcing \MuSeq is defined coordinate wise, i.e.,
for $p,q \in \Forcing \MuSeq$, $q \leq p$ if $\Dom p \Subset \Dom q$
and for every $\xi \in \Dom p$, $q(\xi) \leq p(\xi)$.

\end{itemize}
The weakest condition in \Forcing \MuSeq is the empty function,
denoted by \One. For each $p \in \Forcing \MuSeq$ and $\xi \in \Dom
p$, we let the condition $p(\xi)$ be the pair \Pair {\Tree p \xi}
{\Seq {\br p \xi \delta} {\delta \in \Ps p \xi}}. From now on, \Ps p -
denotes the set \Set {\Pair \xi \delta} {\xi \in \Dom p \And \delta
\in \Ps p \xi}.
\end{DEFINITION}

\begin{FACT}{prod_basic}%
\begin{ITEMS}%
\ITEM{cc} The forcing \Forcing \MuSeq is $\kappa$-closed and it has \SuccCard
\kappa-c.c. 

\ITEM{no_new_branches}%
Suppose $G$ is a \Forcing \MuSeq-generic set over $V$. In \GenExt G,
for every $\xi < \lambda$, the $\kappa$-branches through the tree
$\Tree G \xi = \BigUnion {p \in G} {\Tree p \xi}$ are \Set {\br G \xi
\delta} {\delta < \mu_\xi}, where each $\br G \xi \delta$ is the
function
\[
	\BigUnion - {\Set {\br p \xi \delta}
	{p \in G, \xi \in \Dom p \And
	\delta \in \Ps p \xi}}.
\]
\end{ITEMS}
\begin{Proof}%
\ProofOfItem {no_new_branches}%
Since $\One \Forces {\Res \MuSeq {(\xi+1)}} \Quote {\Forcing {\Res
\MuSeq {(\kappa \Minus (\xi+1))}} \Text {is $\kappa$-closed}}$, the
claim follows from \Lemma{no_new_branches}.%
\end{Proof}%
\end{FACT}

\begin{DEFINITION}{Pairs}%
\Lapikayty{8.10.}%
For all \Forcing \MuSeq-names $\tau$, define that
\[
	\Ps \tau - = \BigUnion - {\Set {\Ps p -}
	{\Text{condition $p$ appears in $\tau$}}}.
\]
Let \Ps \tau 1 denote the set \Set {\xi} {\Pair \xi \delta \in \Ps
\tau -} and \Ps \tau \xi denote the set \Set {\delta} {\Pair \xi
\delta \in \Ps \tau -}.
\end{DEFINITION}

\begin{DEFINITION}{Subforcing}%
\Lapikayty{1.10.}%
Suppose $\ZSeq = \Seq {z_\xi} {\xi \in Z}$ is a sequence such that $Z
\Subset \lambda$ and for each $\xi \in Z$, $z_\xi$ is a subset of
$\mu_\xi$ of cardinality at least $\kappa$. In order to keep our
notation coherent, let \Ps \ZSeq - be a shorthand for the set
\BigUnion {\xi \in Z} {\Braces \xi \times z_\xi}.  We define
\[
	\Forcing \ZSeq = \Set {p \in \Forcing \MuSeq}
	{\Ps p - \Subset \Ps \ZSeq -}.
\]
We say that \Forcing \ZSeq is a subforcing of \Forcing \MuSeq when
\ZSeq is a sequence as described above.
\end{DEFINITION}

A forcing $Q$ is \Def {a complete subforcing of $P$} if every maximal
antichain in $Q$ is also a maximal antichain in $P$ \Note {a set $X$
of conditions is an antichain in $Y$ if all $p \not= q$ in $X$ are
incompatible, i.e., there is no $r \in Y$ with $r \leq p,q$}. The
following basic facts we need later on.

\begin{FACT}{Complete_subforcing}%
\begin{ITEMS}%
\ITEM{complete}%
Every subforcing \Forcing \ZSeq is a complete subforcing of \Forcing
\MuSeq.

\ITEM{one_forces}%
For every $p \in \Forcing \MuSeq$, the restriction \Set {q \in
\Forcing \MuSeq} {q \leq p} is a forcing notion which is equivalent to
\Forcing \MuSeq.

\end{ITEMS}%
\end{FACT}


The following two definitions will be our main tools. Namely, every
permutation $\pi$ of the indices of the labels of the branches in the
generic trees added by \Forcing \MuSeq determines an automorphism
\Autom \pi of \Forcing \MuSeq. This means that for every condition $p$
in \Forcing \MuSeq and \Forcing \MuSeq-name $\tau$ there are many
``isomorphic'' copies of $p$ and $\tau$ inside \Forcing
\MuSeq. Naturally, the copies $\Autom \pi (p)$ and $\Autom \pi (\tau)$
of $p$ and $\tau$, respectively, satisfies all the same formulas \Note
{see \Property {isomorphic} below}.

\begin{DEFINITION}{Maps}%
\Lapikayty{8.10.}%
We define \Maps \MuSeq to be the set of all functions $\pi$ which can
be defined as follows. The domain of $\pi$ is \Ps \YSeq - for some
sequence $\YSeq = \Seq {y_\xi} {\xi \in Y}$ with $Y \Subset \lambda$
and $y_\xi \Subset \mu_\xi$ for each $\xi \in Y$. In addition, there
exists an injective function \Map \pi 1 - from $Y$ into $\lambda$ and
injective functions \Map \pi \xi - from $y_\xi$ into $\mu_\xi$, for
all $\xi \in Y$, such that for all $\Pair \xi \delta \in \Dom \pi$,
\[
	\pi \Pair \xi \delta =
	\Pair {\Map \pi 1 \xi} {\Map \pi \xi \delta}.
\]
When \Forcing \ZSeq is a subforcing of \Forcing \MuSeq, let \Maps
\ZSeq be the collection \Set {\pi \in \Maps \MuSeq} {\Dom \pi \Subset
\Ps \ZSeq -}.
\end{DEFINITION}

\begin{DEFINITION}{images_of_cond_and_names}%
\Lapikayty{8.10.}%
For every $p \in \Forcing \MuSeq$ and $\pi \in \Maps \MuSeq$ with $\Ps
p - \Subset \Dom \pi$, we let $\pi(p)$ denote the condition $q$ in
\Forcing \MuSeq for which
\begin{itemize}

\item $\Dom q = \Map \pi 1 - [\Dom p]$,

\item for every $\zeta \in \Dom q$, $\Tree q \zeta = \Tree p \xi$ and
$\Ps q \zeta = \Map \pi \xi - [\Ps p \xi]$, where $\xi = \Inv {(\Map
\pi 1 -)} (\zeta)$;

\item for every $\Pair \zeta \epsilon \in \Ps q -$, $\br q \zeta
\epsilon = \br p \xi \delta$, where $\Pair \xi \delta = \Inv \pi \Pair
\zeta \epsilon$.

\end{itemize}
When $\tau$ is a \Forcing \ZSeq-name and $\pi$ a mapping in \Maps
\MuSeq with $\Ps \tau - \Subset \Dom \pi$, $\pi(\tau)$ denotes the
\Forcing \ZSeq-name which is result of recursively replacing every
condition $p$ in $\tau$ with $\pi(p)$, i.e.,
\[
	\pi(\tau) =
	\Set {\Pair {\pi(\sigma)} {\pi(p)}} {\Pair \sigma p \in \tau}.
\]
Analogously, for sequences $\ZSeq = \Seq {z_\xi} {\xi \in Z}$ with
$\Ps \ZSeq - \Subset \Dom \pi$, we let $\pi(\ZSeq)$ denote the
sequence \Seq {z_\zeta'} {\zeta \in Z'}, where $Z' = \Map \pi 1 - [Z]$
and for each $\zeta \in Z'$, $z_\zeta' = \Map \pi \xi - [z_\xi]$ with
$\xi = \Inv {(\Map \pi 1 -)} (\zeta)$.
\end{DEFINITION}

\begin{FACT}{isomorphism}%
\Lapikayty{8.10.}%
For every subforcing \Forcing \ZSeq and $\pi \in \Maps \ZSeq$ with
$\Dom \pi = \Ps \ZSeq -$, the mapping \Mapping p {\pi(p)} is an
isomorphism between \Forcing {\ZSeq} and \Forcing {\pi(\ZSeq)}.
\end{FACT}

Suppose \Forcing \ZSeq is a subforcing of \Forcing \MuSeq. The
isomorphism determined by some $\pi \in \Maps \ZSeq$ is denoted by
\Autom \pi. It follows that if $\Dom \pi = \Ps \ZSeq -$, $p \in
\Forcing \ZSeq$, $\psi(x_1, \dots, x_n)$, $n < \omega$, is any
formula, and $\tau_1$, \dots, $\tau_n$ are \Forcing \ZSeq-names then
\PROPERTY{isomorphic}{%
	p \Forces \ZSeq \psi(\tau_1, \dots, \tau_n)
	\Iff
	\Autom \pi(p) \Forces {\pi(\ZSeq)}
	\psi(\Autom \pi(\tau_1), \dots, \Autom \pi(\tau_n)).
}%
Particularly, a mapping $\pi$ in \Maps \ZSeq determines an
automorphism of \Forcing \ZSeq when \Map \pi 1 - is a permutation of
$Z$ and each \Map \pi \xi - is a bijection from $z_\xi$ onto $z_{\Map
\pi 1 \xi}$.
%

\end{SECTION}

\begin{SECTION} {-} {Sigma11}{\protect %
		Basic facts on $\Sigma_1^1$-equivalence relations
}
%
%
Recall that we assumed $\kappa$ to be an uncountable regular cardinal.
We denote the set of all sets hereditarily of cardinality $< \kappa$
by \Her \kappa, i.e., \Her \kappa contains all the sets whose
transitive closure has cardinality $< \kappa$.

\begin{DEFINITION}{Sigma}%
We say that $\phi$ defines a $\Sigma^1_1$-equivalence relation \Equiv
on \Functions \kappa 2 with a parameter $\Param \Subset \Her \kappa$
when%
\begin{ITEMS}%
\ITEM{voc}%
$\phi$ is a first order sentence in the vocabulary consisting of $\in$,
one unary relation symbol $S_0$, and binary relation symbols $S_1$,
$S_2$, and $S_3$;

\ITEM{eq}%
the following definition gives an equivalence relation on \Functions
\kappa 2: for all $f,g \in \Functions \kappa 2$ %
 \[
	f \Equiv [\phi,\Param] g \Iff \ 
	\ForSome h \in \Functions \kappa 2,
	\HerStr \kappa {\Param, f, g, h} \models \phi,
 \]
where \Param, $f$, $g$, and $h$ are the interpretations of the symbols
$S_0$, $S_1$, $S_2$, and $S_3$ respectively.

\end{ITEMS}
We abbreviate \Card [\big] {\Set {\Quotient f {\Equiv}} {f \in
\Functions \kappa 2}} by \No \Equiv.
\end{DEFINITION}

\begin{LEMMA}{Tree_and_Eq}%
\begin{ITEMS}%
\ITEM{less_than_kappa}%
For every nonzero cardinal $\mu \in \kappa \Union \Braces {\kappa,
2^\kappa}$, there exists a $\Sigma_1^1$-equivalence relation \Equiv on
\Functions \kappa 2 with $\No \Equiv = \mu$.

\ITEM{kappa+}%
There exists a $\Sigma_1^1$-equivalence relation \Equiv on \Functions
\kappa 2 with $\No \Equiv = \SuccCard \kappa$.

\ITEM{tree}%
If $T$ is a tree with $\Card T = \kappa$, then there exists a
$\Sigma^1_1$-equivalence relation \Equiv on \Functions \kappa 2 with
$\No \Equiv = \Card {\Br \kappa T} + 1$.

\end{ITEMS}

\begin{Proof}%
Let $\rho$ be a fixed definable bijection from $\kappa$ onto $\kappa
\times \kappa$. For a binary relation $R$, we denote the set \Set
{\rho(\xi)} {\ForSome \xi < \kappa, \Pair \xi 1 \in R} by $\rho(R)$.

\ProofOfItem{less_than_kappa}%
In the cases $\mu \in \kappa \Union \Braces \kappa$, the parameter can
code a list of $\mu$-many nonequivalent functions. In the case $\No
\Equiv = 2^\kappa$ all the functions in \Functions \kappa 2 can be
nonequivalent.

\ProofOfItem{kappa+}%
A sentence $\phi(R_1, R_2, R_3)$ saying
\begin{quote}%
``(both $\rho(R_1)$ and $\rho(R_2)$ are well-orderings of $\kappa$,
and $\rho(R_3)$ is an isomorphism between them) or (neither
$\rho(R_1)$ nor $\rho(R_2)$ is a well-ordering of $\kappa$)''
\end{quote}%
defines a $\Sigma^1_1$-equivalence relation as wanted.

\ProofOfItem{tree}%
We may assume, without loss of generality, that the elements of $T$
are ordinals below $\kappa$. Using \Structure {T, <} as a parameter,
let a sentence $\phi(R_0,R_1, R_2)$ \Note {see \Definition {Sigma}}
say that
\begin{quote}%
``($\rho(R_1) = \rho(R_2)$ is a $\kappa$-branch in $R_0$) or (neither
$\rho(R_1)$ nor $\rho(R_2)$ is a $\kappa$-branch in $R_0$)''.
\end{quote}%
Then $\phi$ defines a $\Sigma^1_1$-equivalence relation as wanted.
\end{Proof}
\end{LEMMA}

\begin{CONCLUSION}{existence}
\Lapikayty{8.10.}%
Let $G$ be a \Forcing \MuSeq-generic set over $V$. Then in \GenExt G,
for every nonzero cardinal $\chi$ in $\kappa \Union \Braces {\kappa,
\SuccCard \kappa, 2^\kappa} \Union \Set {\mu_\xi} {\xi < \lambda}$,
there exists a $\Sigma_1^1$-equivalence relation \Equiv with $\No
\Equiv = \chi$.

\begin{Proof}
The claim follows from \Fact {prod_basic} together with \Lemma
{Tree_and_Eq}.
\end{Proof}
\end{CONCLUSION}

In the next section we shall need the following properties of
$\Sigma_1^1$-equivalence relations.


\begin{LEMMA}{Eq_closure}%
Suppose $\gamma \leq \kappa$ and $\chi_i$, $i < \gamma$, are nonzero
cardinals such that $\phi_i$ defines a $\Sigma_1^1$-equivalence
relation on \Functions \kappa 2 with the parameter $R_i$ and it has
exactly $\chi_i$-many equivalence classes.%
\begin{ITEMS}%
\ITEM{union}%
There exists a $\Sigma_1^1$-equivalence relation \Equiv on \Functions
\kappa 2 with $\No \Equiv = \BigUnion {i < \gamma} {\chi_i}$.

\ITEM{product}%
There exists a $\Sigma_1^1$-equivalence relation \Equiv on \Functions
\kappa 2 with $\No \Equiv = \Card {\Product {i < \gamma} {\chi_i}}$.

\end{ITEMS}

\begin{Proof}%
Both of the claims are simple corollaries of the fact that there are a
parameter $\Param \Subset \Her \kappa$ and a formula $\psi(x)$ such
that for all $f,g,h \in \Functions \kappa 2$
\[
	\HerStr \kappa {\Param, f, g, h}
	\models \psi(i)
\]	
if, and only if,
\[
	\HerStr \kappa {\Param[i], f[i], g[i], h[i]}
	\models \phi_i,
\]	
where $\Param[i]$, $f[i]$, $g[i]$, and $h[i]$ are the \Th i parts of
\Param, $f$, $g$, and $h$ respectively, in some definable coding.
Furthermore $\Param[i] = R_i$ holds for every $i < \gamma$.
\end{Proof}
\end{LEMMA}


\begin{LEMMA}{Absolute}%
Suppose that \Forcing \ZSeq is a subforcing of \Forcing \MuSeq, $\phi$
is a sentence as in \ItemOfDefinition {Sigma} {voc}, and \ParamName,
$\sigma_1$, $\sigma_2$ are \Forcing \ZSeq-names of cardinality
$\kappa$ for subsets of \Her \kappa.
\begin{ITEMS}%
\ITEM{Sigma^1_1}%
If $p \in \Forcing \ZSeq$ and $\psi_1$ denotes the sentence
\[
\Quote{
	\ThereIs h \in \Functions \kappa 2 \With
	\HerStr \kappa {\ParamName, \sigma_1, \sigma_2, h}
	\models \phi
},
\]
then $p \Forces \ZSeq \psi_1$ implies that $p \Forces \MuSeq
\psi_1$ holds, too. 

\ITEM{Auxiliary1}%
\Note{An auxiliary fact only applied in \Item {Pi^1_1} of this
lemma.} %
Suppose $\tau$ is a \Forcing \MuSeq-name of cardinality $\kappa$ for a
subset of \Her \kappa and $q$ is a condition in \Forcing \MuSeq
forcing that
\[
	\HerStr \kappa {\ParamName, \sigma_1, \sigma_2, \tau}
	\models \phi.
\]
For any injection $\eta$ from $\Ps \tau 1 \Minus Z$ into $\lambda
\Minus Z$ there is $\rho \in \Maps \MuSeq$ satisfying that $\Ps \ZSeq
- \Union \Ps q - \Union \Ps \tau - \Subset \Dom \rho$, \Res \rho {\Ps
\ZSeq -} is identity, $\eta \Subset \Map \rho 1 -$, and
\[
	\rho(q) \Forces \MuSeq
	\HerStr \kappa {\ParamName, \sigma_1, \sigma_2,
	\rho(\tau)} \models \phi.
\]
\Remark The conclusion holds even thought $\rho$ does not determine an
automorphism of \Forcing \MuSeq.

\ITEM{Pi^1_1}%
Suppose the length of \ZSeq is at least \SuccCard \kappa, the
cardinality of each $z_\xi$ is at least $\SuccCard \kappa$, $p \in
\Forcing \ZSeq$, and $\psi_2$ denotes the sentence
\[
\Quote{
	\ForAll h \in \Functions \kappa 2,
	\HerStr \kappa {\ParamName, \sigma_1, \sigma_2, h}
	\models \phi
}.
\]
Then $p \Forces \ZSeq \psi_2$ implies $p \Forces \MuSeq \psi_2$.

\ITEM{Auxiliary2}%
\Note{An auxiliary fact only applied in \Lemma {Subforcing}.} %
Suppose $\tau$ is a \Forcing \MuSeq-name of cardinality $\kappa$ for a
subset of \Her \kappa and $q$ is a condition in \Forcing \MuSeq
forcing that
\[
	\Quote{\ForAll h \in \Functions \kappa 2, 
	\HerStr \kappa {\ParamName, \sigma_1, \tau, h}
	\models \phi
	}.
\]
Then for any injection $\rho$ from $\Ps \tau 1 \Minus Z$ into $\lambda
\Minus Z$ there is $\pi \in \Maps \MuSeq$ satisfying that $\Ps \ZSeq -
\Union \Ps q - \Union \Ps \tau - \Subset \Dom \rho$, \Res \pi {\Ps
\ZSeq -} is identity, $\rho \Subset \Map \pi 1 -$, and
\[
	\pi(q) \Forces \MuSeq \Quote{
		\ForAll h \in \Functions \kappa 2,
		\HerStr \kappa {\ParamName, \sigma_1, \pi(\tau),
		h} \models \phi
	}.
\]

\end{ITEMS}

%
\begin{Proof}%
\ProofOfItem{Sigma^1_1}%
%
The claim follows from the facts that \Forcing \MuSeq does not add new
elements into \Her \kappa, and the truth in \Her \kappa is absolute.

\ProofOfItem{Auxiliary1}%
%
Define the subforcing \Forcing \YSeq to be ``the smallest one''
containing \ZSeq, $q$, and $\tau$, i.e., define \YSeq to be the
sequence $\Seq {y_\xi} {\xi \in Y}$ satisfying $\Ps \YSeq - = \Ps
\ZSeq - \Union \Ps q - \Union \Ps \tau -$. Since the truth in \Her
\kappa is absolute,
\[
	q \Forces \YSeq
	\HerStr \kappa {\ParamName, \sigma_1, \sigma_2, \tau}
	\models \phi.
\]
Now each $y_\xi$, $\xi \not\in Z$, is so small that there is $\rho \in
\Maps \MuSeq$ satisfying the demands:
$\Dom \rho = \Ps \YSeq -$,
\Res {\Map \rho 1 -} Z is identity,
\Res {\Map \rho 1 -} {(\Ps \tau 1 \Minus Z)} is $\eta$,
and for every $\xi \in Y$,
\Map \rho \xi - is identity if $\xi \in Z$,
and otherwise, \Map \rho \xi - is some injection from $y_\xi$ into
$\mu_{\eta(\xi)}$.
Since $\rho$ determines an isomorphism between \Forcing \YSeq and
\Forcing {\rho(\YSeq)}, we have that
\[
	\rho(q) \Forces {\rho(\YSeq)}
	\HerStr \kappa {\ParamName, \sigma_1, \sigma_2,
	\rho(\tau)} \models \phi.
\]
Again, by the absoluteness of the truth in \Her \kappa, we can
conclude that the condition $\rho(q)$ forces the same sentence in the
larger forcing \Forcing \MuSeq.

\ProofOfItem{Pi^1_1}%
%
Assume, contrary to the claim, that $p$ is a condition in \Forcing
\ZSeq forcing $\psi_2$, $q \leq p$ is a condition in \Forcing \MuSeq,
and $\tau$ is a \Forcing \MuSeq-name for a function from $\kappa$ into
2 so that
\EQUATION{not_in_ext}{%
%
	q \Forces \MuSeq \Par [\Big] {
	\HerStr \kappa {\ParamName,
	\sigma_1, \sigma_2, \tau} \not\models \phi
	}.
}%
%

We shall define a mapping $\pi \in \Maps \MuSeq$ such that it
determines an automorphism \Autom \pi of \Forcing \MuSeq with the
following properties: $\Autom \pi (\ParamName) = \ParamName$, $\Autom
\pi (\sigma_1) = \sigma_1$, $\Autom \pi (\sigma_2) = \sigma_2$,
$\Autom \pi (\tau)$ is a \Forcing \ZSeq-name, $\Autom \pi (q) \in
\Forcing \ZSeq$, and $\Autom \pi (q)$ is compatible with $q$. It
follows that
$%
	\Autom \pi (q) \Forces \MuSeq \Par [\Big] {
	\HerStr \kappa {\ParamName,
	\sigma_1, \sigma_2, \Autom \pi (\tau)}
	\not\models \phi
	}.
$ %
Then, by the absoluteness of truth in \Her \kappa, we have that
\[
	\Autom \pi (q) \Forces \ZSeq \Par [\Big] {
	\HerStr \kappa {\ParamName,
	\sigma_1, \sigma_2, \Autom \pi (\tau)}
	\not\models \phi
	}.
\]
Since there exist $r \in \Forcing \MuSeq$ with $r \leq p$ and $r \leq
\Autom \pi (q)$, we have a contradiction.

We need the demands that $\Card Z > \kappa$ and $\Card {z_\xi} >
\kappa$, $\xi \in Z$, to ensure that $Z \Minus (\Ps {\sigma_1} -
\Union \Ps {\sigma_2} -)$ and each $z_\xi \Minus (\Ps {\sigma_1} \xi
\Union \Ps {\sigma_2} \xi)$ has cardinality $\geq \kappa$ \Note
{otherwise it is difficult to choose $\pi$ satisfying both $\pi(\tau)
\Subset \Ps \ZSeq -$ and \Res \pi {(\Ps {\sigma_1} - \Union \Ps
{\sigma_2} -)} is identity}.

\Remark A mapping $\pi \in \Maps \MuSeq$ determines an automorphism if
\Map \pi 1 - is a permutation of $\lambda$ and each \Map \pi \xi - is
a bijection from $\mu_\xi$ onto $\mu_{\Map \pi 1 \xi}$. Thus we need
that the chosen $\pi$ satisfies $\mu_{\Map \pi 1 \xi} = \mu_\xi$ for
every $\xi \in \Ps \tau 1$. Now $Z$ might be too small to contain all
the possible cardinals in \MuSeq. However, because of \Item
{Auxiliary1} and the assumption that all the cardinals in \MuSeq are
repeated $\lambda$-many times, there are $q'$ and $\tau'$ satisfying
\Equation {not_in_ext}. Moreover, 
for every $\xi \in \Ps \tau 1 \Minus Z$, there is $\zeta \in Z \Minus
(\Ps q 1 \Union \Ps \ParamName 1 \Union \Ps {\sigma_1} 1 \Union \Ps
{\sigma_2} 1)$ with $\mu_\zeta = \mu_\xi$, and $q' \leq p$ holds, too
\Note {since $q' = \rho(q)$ and the branches in $p$ are kept fixed,
i.e., \Res \rho {\Ps \ZSeq -} is identity}.

Define, for every $\Pair \xi \delta \in \Ps q - \Union \Ps \tau -$, a
pair \Pair {\zeta_\xi} {\epsilon_{\xi, \delta}} as follows. Set
$\zeta_\xi = \xi$
if $\xi \in \Ps \ParamName 1 \Union \Ps {\sigma_1} 1 \Union \Ps
{\sigma_2} 1$,
and choose
$\zeta_\xi \in Z \Minus (\Ps q 1 \Union \Ps \ParamName 1 \Union \Ps
{\sigma_1} 1 \Union \Ps {\sigma_2} 1)$ with $\mu_{\zeta_\xi} =
\mu_\xi$
otherwise.  Analogously, let $\epsilon_{\xi, \delta}$ be $\delta$ if
$\delta$ is in
$\Ps \ParamName \xi \Union \Ps {\sigma_1} \xi \Union \Ps {\sigma_2}
\xi$,
and pick some $\epsilon_{\xi, \delta}$ from
$z_\xi \Minus (\Ps q \xi \Union \Ps \ParamName \xi \Union \Ps
{\sigma_1} \xi \Union \Ps {\sigma_2} \xi)$
otherwise. Let $\pi$ be any mapping from \Maps \MuSeq which satisfies
that
\begin{itemize}%

\item \Res \pi {(\Ps \ParamName - \Union \Ps {\sigma_1} - \Union \Ps
{\sigma_2} -)} is identity;

\item for every $\Pair \xi \delta \in \Ps q - \Union \Ps \tau -$, $\pi
\Pair \xi \delta = \Pair {\zeta_\xi} {\epsilon_{\xi,\delta}}$.

\item \Map \pi 1 - is a permutation of $\lambda$;

\item for every $\xi < \lambda$, $\Map \pi \xi -$ is a permutation of
$\mu_\xi$;

\end{itemize}%
Then $\pi$ determines an automorphism as wanted.

\ProofOfItem{Auxiliary2}%
%
The proof is similar to the proof of \Item {Auxiliary1}. The main
difference is that one has to apply \Item {Sigma^1_1} and \Item
{Pi^1_1} instead of ``the absoluteness of truth in \Her \kappa''.
\end{Proof}%
%
%
\end{LEMMA}

\begin{CONCLUSION}{Absolute}%
Suppose \Forcing \ZSeq is a subforcing of \Forcing \MuSeq and $\psi$
is a sentence in the vocabulary \Braces {\in, S_0, S_1, S_2} which is
a Boolean combination of a sentence containing one second order
existential quantifier and a sentence containing one second order
universal quantifier.  Then for every \Forcing \MuSeq-generic set $H$
over $V$, $G = H \Inter \Forcing \ZSeq$ is \Forcing \ZSeq-generic set
over $V$, $\GenExt G \Subset \GenExt H$, and for all $f,g \in \InExt
{\Functions \kappa 2} G$,
\[
	\OfGenExt {
		\HerStr \kappa {\Param,f,g} \models \psi
	} G \Iff
	\OfGenExt {
		\HerStr \kappa {\Param,f,g} \models \psi
	} H,
\]
where \Param, $f$, and $g$ are interpretations of the symbols $S_0$,
$S_1$, and $S_2$ respectively.
\end{CONCLUSION}
%
%
\end{SECTION}

\begin{SECTION} {-} {Sigma11No}{\protect %
	Possible numbers of equivalence classes
}
%
\begin{DEFINITION}{ClosedCardSet}%
\Lapikayty{1.10.}%
Suppose $\MuSeq = \Seq {\mu_\xi} {\xi < \lambda}$ is a sequence of
cardinals. Define \ClosedCardSet to be the smallest set of cardinals
satisfying that
\begin{itemize}

\item every nonzero cardinal $\leq \SuccCard \kappa$ is in
\ClosedCardSet;

\item $\Set {\mu_\xi} {\xi < \lambda} \Subset \ClosedCardSet$;

\item if $\gamma \leq \kappa$ and $\chi_i$, $i < \gamma$, are
cardinals in \ClosedCardSet, then both \BigUnion {i < \gamma} {\chi_i}
and \Card {\Product {i < \gamma} {\chi_i}} are in \ClosedCardSet.

\end{itemize}
\end{DEFINITION}

We shall now prove that, when the singular cardinal hypothesis holds,
the closure under unions and products as above are, consistent wise,
the only restrictions on the possible numbers of equivalence classes
of $\Sigma^1_1$-equivalence relations.

\begin{THEOREM}{No}%
\Lapikayty{1.10.}%
Suppose that
\begin{itemize}

\item $\kappa$ is an uncountable cardinal with $\kappa^{<\kappa} =
\kappa$ and $2^\kappa = \SuccCard \kappa$;

\item $\lambda > \SuccCard \kappa$ is a cardinal with
$\lambda^\kappa = \lambda$;

\item $\MuSeq = \Seq {\mu_\xi} {\xi < \lambda}$ and \Forcing \MuSeq
are as in \Definition {ProdForc};

\item \ClosedCardSet is as in \Definition {ClosedCardSet};

\item for every $\chi \in \ClosedCardSet$ with $\chi > \SuccCard
\kappa$ and $\gamma < \kappa$, the inequality $\chi^\gamma \leq
\SuccCard \chi$ holds.

\end{itemize}
Then for every \Forcing \MuSeq-generic set $G$, the extension \GenExt
G satisfies that all cardinals and cofinalities are preserved, there
are no new sets of cardinality $< \kappa$, $2^\kappa = \lambda$ and
for all cardinals $\chi$, the following two conditions are equivalent:
\begin{ITEMS}%
\ITEM{InClosedCardSet} $\chi \in \ClosedCardSet$;

\ITEM{NoEquiv}%
a sentence $\phi$ defines a $\Sigma_1^1$-equivalence relation \Equiv
on \Functions \kappa 2 with a parameter $R \Subset \Her \kappa$ and
there are exactly $\chi$ different equivalence classes of \Equiv.
\end{ITEMS}
\end{THEOREM}

\PvComment{%
Since \Forcing \MuSeq has \SuccCard \kappa-c.c. all the unions of
length $\leq \kappa$ remains the same in the ground model and in the
extension. So \ClosedCardSet is the same in the ground model and in
the extension.
}

The rest of this section is devoted to the proof of this theorem.
Because of \Conclusion {existence} and \Lemma {Eq_closure} it remains
to show that
\[
	\One \Forces \MuSeq \Quote {
	\No {\Equiv[\phi,\Name{S}]} \in \ClosedCardSet
	\Text{for all $\Sigma^1_1$-equivalence
	relations \Equiv[\phi,\Name{S}]}
	}.
\]
Suppose $p$ is a condition in \Forcing \MuSeq and $\theta$ a cardinal
such that
\[
	p \Forces\MuSeq \Quote{
	\ThereExists X \Subset \Her \kappa \With
	\No {\Equiv[\phi,X]} = \theta}.
\]
By \ItemOfFact {Complete_subforcing} {one_forces} we have that the
condition \One forces the same formula. Hence by the maximal principle
we may fix a name \ParamName so that $\One \Forces\MuSeq \ParamName
\Subset \Her \kappa$ and
\EQUATION{initial_assumption}{%
	\One \Forces\MuSeq \No\EquivName = \theta.
}
Since \Forcing \MuSeq has \SuccCard \kappa-c.c. and $\Card {\Her
\kappa} = \kappa^{<\kappa} = \kappa$, we may assume that the name
\ParamName has cardinality $\kappa$.



\subsection{Choice of a small subforcing}

Next we want to prove that there is a subforcing \Forcing \ZSeq of
\Forcing \MuSeq such that the cardinality of \Forcing \ZSeq is
$\theta$, there are only \SuccCard \kappa-many coordinates in \Forcing
\ZSeq, and already \Forcing \ZSeq produces $\theta$-many different
equivalent classes of \EquivName.

\begin{LEMMA}{Subforcing}%
Suppose \Forcing \ZSeq is a subforcing of \Forcing \MuSeq such that%
\begin{itemize}%

\item $\ZSeq = \Seq {z_\xi} {\xi \in Z}$;

\item $Z$ is a subset of $\lambda$ satisfying $\Card Z = \SuccCard
\kappa$ and $\Ps \ParamName 1 \Subset Z$;

\item for each $\xi \in \Ps \ParamName 1$, $z_\xi = \mu_\xi$ if
$\mu_\xi \leq \theta$, and otherwise, $z_\xi \in \Set {y \in \Subsets
{\mu_\xi} {\SuccCard \kappa}} {\Ps \ParamName \xi \Subset y}$.

\item $Z \Minus \Ps \ParamName 1$ is of cardinality \SuccCard \kappa;

\item for every $\xi \in Z \Minus \Ps \ParamName 1$, $\mu_\xi >
\theta$ and $z_\xi$ is some set in \Subsets {\mu_\xi} {\SuccCard
\kappa}.

\end{itemize}
Then our assumption \EquationOnPage [\EquationRefCtrl]
{initial_assumption} implies $\One \Forces\ZSeq \No \EquivName =
\theta$.

%
\begin{Proof}%
Let \FunName\ZSeq be a \Forcing \ZSeq-name for the set of all
functions from $\kappa$ into 2, i.e., it satisfies that $\One \Forces
\ZSeq \FunName\ZSeq = \Functions \kappa 2$. We prove that
\EQUATION{subclaim}{%
	\One \Forces\MuSeq \Quote{
		\ForAll f \in \Functions \kappa 2
		\ThereIs g \in \FunName \ZSeq
		\With f \EquivName g
	}.
}%
This suffices since then
\[
	\One \Forces\MuSeq
	\theta
	\leq \Card {\Quotient {\FunName \ZSeq} {\EquivName}}
	\leq \Card {\Quotient {\Functions \kappa 2} {\EquivName}} 
	= \No \EquivName
	= \theta,
\]
and, by \LemmaItem {Absolute} {Sigma^1_1}, \LemmaItem {Absolute}
{Pi^1_1} of \Lemma {Absolute}, we can conclude
\[
	\One \Forces\ZSeq
	\No \EquivName
	= \Card {\Quotient {\FunName \ZSeq} {\EquivName}}
	= \theta.
\]

Now assume, contrary to \Equation {subclaim}, that \EquationOnPage
[\EquationRefCtrl] {initial_assumption} holds, there is a condition
$p$ in \Forcing \MuSeq, and a \Forcing \MuSeq-name $\sigma$ for a
function from $\kappa$ into 2 such that
\EQUATION{existence_of_noneq_function}{%
	p \Forces \MuSeq \Quote {
		\ForAll g \in \FunName \ZSeq,
		\sigma \NotEquivName g\,
	}.%
}%
We may assume that the name $\sigma$ is of cardinality $\kappa$.
Furthermore, by \ItemOfLemma {Absolute} {Auxiliary2} and since each
cardinal in \MuSeq is listed $\lambda$-many times, we may choose $p$
and the name $\sigma$ so that the coordinates appearing in $\sigma$
adds a tree with the same number of $\kappa$-branches as some
coordinate in $Z$ does, i.e., for every $\xi \in \Ps \sigma 1$, there
is $\zeta \in Z$ with $\mu_\zeta = \mu_\xi$. This property will be
essential in the choice of automorphisms \Note {in the same way as the
analogous demand \ItemOfLemma {Absolute} {Auxiliary1} was needed in
the proof of \ItemOfLemma {Absolute} {Pi^1_1}, see the remark in the
middle of that proof}.

Our strategy will be the following.
\begin{ITEMS}%

\ITEM{sigma'}%
We define a name $\sigma'$ so that $\One \Forces \MuSeq \sigma' \in
\FunName \ZSeq$. Hence, by applying \Equation
{existence_of_noneq_function}, we get
\[
	p \Forces \MuSeq \sigma \NotEquivName \sigma'.
\]

\ITEM{maps}%
We define \Forcing \MuSeq-names \Seq {\tau^\gamma} {\gamma < \SuccCard
\theta} for functions from $\kappa$ into 2, and conditions \Seq
{q^\gamma} {\gamma < \SuccCard \theta} in \Forcing \MuSeq.

\ITEM{maps_properties}%
For every $\gamma < \gamma' < \SuccCard \theta$ we define a mapping
\Map [\gamma, \gamma'] \rho -- in \Maps \MuSeq such that \Map [\gamma,
\gamma'] \rho -- determines an automorphism \Map [\gamma, \gamma']
{\Autom \rho} -- of \Forcing \MuSeq with the following properties:
$\Map [\gamma, \gamma'] {\Autom \rho} - \ParamName = \ParamName$,
$\Map [\gamma, \gamma'] {\Autom \rho} - p = q^\gamma$,
$\Map [\gamma, \gamma'] {\Autom \rho} - \sigma = \tau^\gamma$, and
$\Map [\gamma, \gamma'] {\Autom \rho} - {\sigma'} = \tau^{\gamma'}$.
Hence it follows from \Item {sigma'} that
%
%
\[
	q^\gamma \Forces \MuSeq
	\tau^\gamma \NotEquivName \tau^{\gamma'}.%
\]
%

\ITEM{nonequivalent_in_genext}%
Finally, we fix a \Forcing \MuSeq-generic set $G$ over $V$ and, by
applying ``a standard density argument'', we show that for some $B \in
\Subsets {\SuccCard \theta} {\SuccCard \theta}$, all the conditions
$q^\gamma$, $\gamma \in B$, are in the generic set $G$. It follows
from \Item {maps_properties} that in \GenExt G, $\No \Equiv \geq
\SuccCard \theta$ contrary to \EquationOnPage [\EquationRefCtrl]
{initial_assumption}.

\end{ITEMS}%

As can be guessed from the demands on the sequence \ZSeq, there are
three different kind of indices which we have to deal with:
\begin{itemize}%

\item $\Theta_\leq = \Set {\xi \in \Ps \ParamName 1} {\mu_\xi \leq
\theta}$,

\item $\Theta_> = \Set {\xi \in \Ps \ParamName 1} {\mu_\xi > \theta}$,
and

\item $\Theta' = \lambda \Minus \Ps \ParamName 1$.

\end{itemize}%

\Remark Of course we would like to have that $q^\gamma = \Map [\gamma,
\gamma'] \rho - p = p$ for every $\gamma < \gamma' < \SuccCard
\theta$. Unfortunately, that is not possible since it might be the
case that for some $\xi \in \Theta_>$, $\Ps \sigma \xi \Inter \Ps p
\xi \not\Subset z_\xi$ \Note {and we really need later the restriction
$\Card {z_\xi} < \theta$}.

\ProofOfItem{sigma'}%
We define the name $\sigma'$ to be $\Map \pi - \sigma$ for a mapping
$\pi$ in \Maps \MuSeq which satisfies the following conditions:
\begin{itemize}%

\item $\Dom \pi = \Ps \sigma -$;

\item $\Ran \pi \Subset \Ps \ZSeq -$;

\item for every $\xi \in \Dom {\Map \pi 1 -}$, $\mu_{(\Map \pi 1 \xi)}
= \mu_\xi$;

\item \Res \pi {\Ps \ParamName -} is identity \Note {implying \Res
{\Map \pi 1-} {(\Theta_\leq \Union \Theta_>)} is identity};

\item for every $\xi \in \Dom {\Map \pi 1 -} \Inter \Theta_\leq$, \Map
\pi \xi - is identity;

\item for every $\xi \in \Dom {\Map \pi 1 -} \Inter \Theta_>$ and
$\delta \in \Dom {\Map \pi \xi -} \Minus \Ps \ParamName \xi$, $\Map
\pi \xi \delta \not\in \Ps p \xi \Union \Ps \ParamName \xi \Union \Ps
\sigma \xi$;

\item for every $\xi \in \Dom {\Map \pi 1 -} \Inter \Theta'$, $\Map
\pi 1 \xi \not\in \Dom p \Union \Ps \ParamName 1 \Union \Ps \sigma 1$
and $\pi_\xi$ is some injective function having range $z_\xi$.

\end{itemize}%
It is possible to fulfill these conditions by the choice of $\sigma$,
because of the cardinality demands on \ZSeq, and since $\Ps p - \Union
\Ps \ParamName - \Union \Ps \sigma -$ has cardinality $\kappa$. Since
$\One \Forces \MuSeq \sigma \in \Functions \kappa 2$ and $\pi$ can be
extended so that the extension determines an automorphism of \Forcing
\MuSeq, we have that $\One \Forces \MuSeq \sigma' \in \Functions
\kappa 2$. However, $\sigma'$ is a \Forcing \ZSeq-name, so $\One
\Forces \MuSeq \sigma' \in \FunName \ZSeq$ holds, too.

\ProofOfItem{maps}%
For every $\gamma < \SuccCard \theta$, we define a mapping $\Map
[\gamma] \pi-- \in \Maps \MuSeq$ so that the desired name
$\tau^\gamma$ is \Map [\gamma] \pi - \sigma and the condition
$q^\gamma$ is \Map [\gamma] \pi - p. Since we do NOT demand that $\Ran
{\Map [\gamma] \pi --} \Subset \Ps \ZSeq -$, when $\gamma < \SuccCard
\theta$, it is possible to choose \Map [\gamma] \pi -- so that all the
following demands are fulfilled:
\begin{itemize}%

\item $\Dom {\Map [\gamma] \pi--} = \Ps \sigma - \Union \Ps p -$

\item \Res {\Map [\gamma] \pi--} {\Ps \ParamName -} is identity;

\item for every $\xi \in \Dom {\Map [\gamma] \pi 1 -}$, $\mu_{(\Map
[\gamma] \pi 1 \xi)} = \mu_\xi$;

\item for every
$\xi \in \Dom {\Map [\gamma] \pi 1 -} \Inter (\Theta_\leq \Union
\Theta')$,
\Map [\gamma] \pi \xi - is identity;

\item for all
$\xi \in \Dom {\Map [\gamma] \pi 1 -} \Inter \Theta_>$,
the sets
$(\Ps \ParamName \xi \Union \Ps \sigma \xi \Union \Ps {\sigma'} \xi
\Union \Ps p \xi)$
and
$\Ran {\Map [\gamma] \pi \xi -} \Minus \Ps \ParamName \xi$,
for all $\gamma < \SuccCard \theta$, are pairwise disjoint;

\item the sets
$(\Ps \ParamName 1 \Union \Ps \sigma 1 \Union \Ps {\sigma'} 1 \Union
\Dom p)$
and
$\Ran {\Map [\gamma] \pi 1 -} \Minus \Ps \ParamName 1$,
for all $\gamma < \SuccCard \theta$, are pairwise disjoint.

\end{itemize}

\ProofOfItem{maps_properties}%
Fix indices $\gamma < \gamma' < \SuccCard \theta$. Consider the set of
pairs \Pair x y satisfying that
\begin{itemize}%

\item		$x \in \Dom {\Map [\gamma] \pi --}$
		and
		$\Map [\gamma] \pi - x = y$, or

\item there is $z \in \Dom \pi = \Ps \sigma -$ such that
		$\pi(z) = x$
		and
		$\Map [\gamma'] \pi - z = y$.

\end{itemize}%
Because of the conditions given above, we have that
\begin{itemize}%

\item for all $\xi \in \Dom {\Map[\gamma] \pi 1-} = \Dom
{\Map[\gamma'] \pi 1-}$, 
$\Map [\gamma] \pi 1 \xi = \Map [\gamma'] \pi 1 \xi$ iff $\Map \pi 1
\xi = \xi$;

\item for all $\Pair \xi \delta \in \Dom {\Map[\gamma] \pi --} = \Dom
{\Map[\gamma'] \pi --}$, $\Map [\gamma] \pi \xi \delta = \Map
[\gamma'] \pi \xi \delta$ iff $\Map \pi \xi \delta = \delta$;

\item for all $\xi \not = \zeta \in \Dom {\Map[\gamma] \pi 1-}$,
$\Map [\gamma] \pi \xi \delta \not= \Map [\gamma'] \pi \zeta \epsilon$;

\item for all $\Pair \xi \delta \not = \Pair \xi \epsilon \in \Dom
{\Map [\gamma] \pi --}$, $\Map [\gamma] \pi \xi \delta \not= \Map
[\gamma'] \pi \xi \epsilon$.

\end{itemize}%
Hence the set of pairs we considered is the following well-defined
injective function from \Maps \MuSeq:
\[
	\eta =
	\Map [\gamma] \pi --
	\Union
	\Par [\Big] {
		(\Res {\Map [\gamma'] \pi --} {\Dom \pi})
		\Comp
		\Inv {(\pi)}
	}.
\]
We let the mapping \Map [\gamma, \gamma'] \rho -- be any extension of
$\eta$ satisfying that $\Map [\gamma, \gamma'] \rho -- \in \Maps
\MuSeq$, $\Dom {\Map [\gamma, \gamma'] \rho 1-} = \lambda$, and for
each $\xi < \lambda$, $\Dom {\Map [\gamma, \gamma'] \rho \xi-} =
\mu_\xi$. It follows that
\begin{itemize}%

\item $\Map [\gamma, \gamma'] \rho - \ParamName = 
\pi(\ParamName) =
\Map [\gamma] \pi - \ParamName = 
\Map [\gamma'] \pi - \ParamName = \ParamName$;

\item $\Map [\gamma, \gamma'] \rho - p =
\Map [\gamma] \pi - p =
q^\gamma$
\Note{note, that $\Ran \pi \Inter (\Ps p- \Minus \Ps \ParamName-) =
\emptyset$};

\item $\Map [\gamma, \gamma'] \rho - \sigma =
\Map [\gamma] \pi - \sigma =
\tau^\gamma$
\Note{note, that $\Ran \pi \Inter (\Ps \sigma- \Minus \Ps \ParamName-)
= \emptyset$};

\item $\Map [\gamma, \gamma'] \rho - {\sigma'} = 
\Map [\gamma'] \pi -- \Par [\Big] {\Inv \pi (\sigma')} =
\Map [\gamma'] \pi -- (\sigma) =
\tau^{\gamma'}$
\Note {remember, that $\pi(\sigma) = \sigma'$}.

\end{itemize}

\ProofOfItem{nonequivalent_in_genext}%
Our demands on the mappings \Map [\gamma] \pi --, $\gamma < \SuccCard
\theta$, ensure that for each $\Pair \xi \delta \in \Ps p -$, if
$\Pair \xi \delta \in \Ps {(q^\gamma)} -$ then $\br {(q^\gamma)} \xi
\delta = \br p \xi \delta$. Therefore, $p$ and $q^\gamma$ are
compatible conditions%
%
%
. Moreover, for every $\beta < \SuccCard \theta$, the set
\[
	D_\beta = \Set {r \in \Forcing \MuSeq} {%
	\ForSome \gamma > \beta, r \leq q^\gamma}
\]
is a dense set below the condition $p$ \Note {which means that for
every $s \leq p$ there is $r \leq s$ with $r \in D_\beta$}.  Because
of $p \in G$, $D_\beta \Inter G$ is nonempty for every $\beta <
\SuccCard \theta$.  Consequently, the set $B = \Set {\gamma <
\SuccCard \theta} {q^\gamma \in G}$ must be cofinal in \SuccCard
\theta. So $B$ has cardinality \SuccCard \theta.%
\end{Proof}%
%
%
\end{LEMMA}



%
\subsection{Isomorphism classes of names}\label{isom_names}

First of all we fix \ZSeq so that the subforcing \Forcing \ZSeq of
\Forcing \MuSeq satisfies the assumptions of \Lemma {Subforcing}.
Secondly we fix \Forcing \ZSeq-names \Seq {\sigma_\alpha} {\alpha <
\theta} for functions from $\kappa$ into 2 so that for all $\alpha
\not= \beta < \theta$,
\EQUATION{noneq_names}{%
	\One \Forces\ZSeq
	\sigma_\alpha \NotEquivName \sigma_\beta.
}%
Since \Forcing \ZSeq has \SuccCard \kappa-c.c., we may assume that
each of the names $\sigma_\alpha$ has cardinality $\kappa$.

\begin{DEFINITION}{isomorphic_names}%
For every $\alpha < \theta$ we fix an enumeration \Seq {\SEp \alpha i}
{i < \kappa} of \Ps {\sigma_\alpha} - without repetition. Names
$\sigma_\alpha$ and $\sigma_\beta$ are said to be isomorphic, written
$\sigma_\alpha \Isomorphic \sigma_\beta$, if the following conditions
are met:
\begin{itemize}

\item for every $i < \kappa$, $\xi^\alpha_i = \xi^\beta_i$;

\item for every $i < \kappa$ and $\zeta = \xi^\alpha_i = \xi^\beta_i$,
if $\mu_\zeta > \theta$ then also $\delta^\alpha_i = \delta^\beta_i$;

\item for all $\Pair \zeta \epsilon \in \Ps \ParamName -$ and $i <
\kappa$, $\SEp \alpha i = \Pair \zeta \epsilon$ iff $\SEp \beta i =
\Pair \zeta \epsilon$.

\item $\pi(\sigma_\alpha) = \sigma_\beta$ when $\pi \in \Maps \ZSeq$
is the mapping with $\Dom \pi = \Ps {\sigma_\alpha} -$ and $\pi \SEp
\alpha i = \SEp \beta i$ for each $i < \kappa$.

\end{itemize}
\end{DEFINITION}

For every $\alpha < \theta$ we denote the set \Set {\beta < \theta}
{\sigma_\beta \Isomorphic \sigma_\alpha} by \IClass [\alpha]. Now by
the choice of \Forcing \ZSeq, and the assumptions $\kappa^{< \kappa} =
\kappa$ and $2^\kappa = \SuccCard \kappa$, the number of nonisomorphic
names in \Set {\sigma_\alpha} {\alpha < \theta} is $\leq \SuccCard
\kappa$, i.e., the cardinality of the family \Set {\IClass [\alpha]}
{\alpha < \theta} is at most \SuccCard \kappa.

\begin{EasyComment}%
\begin{EasyNote}%
\newcommand{\RInd}[1]{\MathX{%
	J^{#1}%
}}
For each $\alpha < \theta$ there is a sequence $\XiSeq \alpha = \Seq
{\xi^\alpha_i} {i < \kappa}$ such that for every $\beta \in \IClass
[\alpha]$, $\XiSeq \beta = \XiSeq \alpha$. The number of different
sequences \XiSeq \alpha is the number of functions from $\kappa$ into
$Z$. Since $Z$ is chosen to have cardinality \SuccCard \kappa, the
number of different such sequences is ${(\SuccCard \kappa)}^\kappa =
\SuccCard \kappa$.

For each $\alpha < \theta$ there is a set $\RInd \alpha \Subset
\kappa$ such that for every $\beta \in \IClass [\alpha]$ and $i <
\kappa$, $i \in \RInd \alpha$ iff $\mu_{\xi^\beta_i} > \theta$ or
$\SEp \beta i \in \Ps \ParamName -$. The number of different sets
\RInd \alpha is $2^\kappa = \SuccCard \kappa$.

For each $\alpha < \theta$ there is a sequence $\EpsilonSeq \alpha =
\Seq {\epsilon^\alpha_i} {i \in \RInd \alpha}$ such that for every
$\beta \in \IClass [\alpha]$, $\Res {\DeltaSeq \beta} {\RInd \alpha} =
\EpsilonSeq \alpha$.  The number of different sequences \EpsilonSeq
\alpha is the number of functions from \RInd \alpha into \Ps
\ParamName -, i.e., $\kappa^\kappa = \SuccCard \kappa$.

Let \Seq {\Pair {\zeta_i} {\epsilon_i}} {i < \kappa} be fixed
enumeration of $\kappa \times \kappa$. Define $\rho^\alpha$ to be the
mapping \Mapping {\Pair {\xi^\alpha_i} {\delta^\alpha_i}} {\Pair
{\zeta_i} {\epsilon_i}}. Then every name $\rho^\alpha(\sigma_\alpha)$,
abbreviated by $\tau^\alpha$, is a \Forcing \YSeq-name, where $\YSeq =
\Seq {y_\xi} {\xi < \kappa}$ is such that every $y_\xi$ is
$\kappa$. By the assumption $\kappa^{<\kappa} = \kappa$, \Forcing
\YSeq has cardinality $\kappa$, and consequently, there are at most
$\SuccCard \kappa$-many different names in \Set {\tau^\alpha} {\alpha
< \theta}.

Then for all $\alpha, \beta < \theta$, $\sigma_\alpha \Isomorphic
\sigma_\beta$ iff $\RInd \alpha = \RInd \beta$, $\XiSeq \alpha =
\XiSeq \beta$, $\EpsilonSeq \alpha = \EpsilonSeq \beta$, $\tau^\alpha
= \tau^\beta$.  Therefore, the number of nonisomorphic names in \Set
{\sigma_\alpha} {\alpha < \theta} is the number of different tuples of
the form \SimpleSeq {\XiSeq \alpha, \RInd \alpha, \EpsilonSeq \alpha,
\tau^\alpha}, which is \SuccCard \kappa.
\end{EasyNote}
\end{EasyComment}

Let \Repr be a subset of $\theta$ such that $\Card \Repr \leq
\SuccCard \kappa$ and \Set {\sigma_\alpha} {\alpha \in \Repr} is a set
of representatives of the isomorphism classes. If $\theta \leq
\SuccCard \kappa$ then $\theta \in \ClosedCardSet$ directly by the
definition. From now on we assume that $\theta > \SuccCard
\kappa$. Hence $\theta = \BigUnion {\alpha \in \Gamma} {\IClass
[\alpha]}$ implies that
\EQUATION{sup_of_isom_classes}{%
	\theta = \BigUnion {\alpha \in \Repr}
	{\Card {\IClass[\alpha]}}.
}%
Define ``the set of all small cardinals'' to be
\[
	\SmallCards =
	\Set {\mu_\xi}
	{\xi \in \Ps \ParamName 1 \And \mu_\xi \leq \theta}.
\]
Note that this set might be empty. Anyway, then we know that
\EQUATION{lower_bound}{%
	\theta \geq \Max {\SuccCard [2] \kappa, \sup \SmallCards}.
}%
So to prove that $\theta$ is a cardinal in \ClosedCardSet we shall
show that for every $\alpha \in \Repr$, the cardinality of \IClass
[\alpha] is strictly smaller than the lower bound given in \Equation
{lower_bound} above, or otherwise, we can find a subset
\CritInd[\alpha] of $\kappa$ so that \Card {\IClass [\alpha]} has one
of the following form:
either $\Card {\CritInd [\alpha]} < \kappa$ and
\EQUATION{small_CritInd}{%
	\Card {\IClass [\alpha]} \in
	\Braces[\Big]{
		\BigUnion {i \in \CritInd[\alpha]} {\mu_{\xi^\alpha_i}}
	}
	\Union
	\Braces[\Big]{
		\Product {i \in \CritInd[\alpha]} {\mu_{\xi^\alpha_i}}
	},%
}%
or else, $\Card {\CritInd[\alpha]} = \kappa$ and 
\EQUATION{large_CritInd}{%
	\Card {\IClass [\alpha]} =
	\BigUnion {K \in \Subsets {\CritInd[\alpha]} {<\kappa}}
	{\Card[\Big]{
		\Product{i \in K}{\mu_{\xi^\alpha_i}}
	}}.
}%
This will suffice since we take care of that for every $\alpha \in
\Repr$ and for each $i \in \CritInd [\alpha]$, $\mu_{\xi^\alpha_i} \in
\SmallCards$, i.e., only those small cardinals are used whose
coordinate appears in the name \ParamName. Then there occurs at most
$\kappa$-many different cardinals in the union \Equation
{sup_of_isom_classes}, and hence, for some sequence \Seq {X_k} {k <
\kappa} of sets in \Subsets {\SmallCards} {<\kappa},
\[
	\theta = \BigUnion {k < \kappa}
		{\Card [\Big] {\Product {\mu \in X_k} {\mu}}}
	\in \ClosedCardSet.
\]

\Remark From our assumption that for every $\chi \in (\ClosedCardSet
\Minus \SuccCard [2] \kappa)$ and $\gamma < \kappa$, the inequality
$\chi^\gamma \leq \SuccCard \chi$ holds, it follows that $\theta$ is
either $\sup \SmallCards$ or \Card [\Big] {\Product {\mu \in X} {\mu}}
for some subset $X$ of \SmallCards with $\Card X < \kappa$.

\subsection{\protect %
Case 1: The parameter depends on $< \kappa$-many coordinates}%
\label{small}%

For the rest of the proof, let $\alpha^*$ be a fixed ordinal so that
the number of names in \Set {\sigma_\beta} {\beta < \theta}, which are
isomorphic to the representative $\sigma_{\alpha^*}$, is greater or
equal to the lower bound given in \EquationOnPage {lower_bound}, i.e.,
$\alpha^* \in \Gamma$ and \Card {\IClass [\alpha^*]} is large enough.
To simplify our notation, let $\XiSeq* = \Seq {\xi^*_i} {i < \kappa}$
and $\DeltaSeq* = \Seq {\delta^*_i} {i < \kappa}$ denote the sequences
\XiSeq {\alpha^*} and \DeltaSeq {\alpha^*} respectively, and
abbreviate \IClass [\alpha^*] by \IClass.

Define the set of ``all critical indices of the isomorphism class of
$\sigma_{\alpha^*}$'' to be
\EQUATION{CritJ}{%
	\CritJ = \Set {i < \kappa}{
	\mu_{\xi^*_i} \leq \theta
	\And 
	\Pair {\xi^*_i} {\delta^*_i}
	\not\in \Ps \ParamName -
	}.
}%
Note that for every $\alpha \in \IClass$, the equations $\XiSeq \alpha
= \XiSeq*$ and $\EqualRes {\DeltaSeq \alpha} {\DeltaSeq*} {(\kappa
\Minus \CritJ)}$ hold. Note also, that by the choice of \Forcing
\ZSeq, $\CritJ \Subset \Set {i < \kappa} {z_{\xi^*_i} = \mu_{\xi^*_i}}
\Subset \Set {i < \kappa} {\xi^*_i \in \Ps \ParamName 1}$. Thus $\Set
{\mu_{\xi^*_i}} {i \in \CritJ} \Subset \SmallCards$ holds, too.

The set \CritJ must be nonempty, since otherwise there are $\alpha
\not= \beta$ in \IClass such that $\sigma_\alpha$ is the same name as
$\sigma_\beta$, contrary to the choice that $\sigma_\alpha$ and
$\sigma_\beta$ are names for nonequivalent functions \Note
{\EquationOnPage {noneq_names}}. For a similar reason $\Card [\big]
{\Product {i \in \CritJ} {\mu_{\xi^*_i}}} \geq \Card \IClass$ holds.

Now suppose that already some subset $K$ of \CritJ having cardinality
$< \kappa$ satisfies the following inequality:
\[
	\Card [\big] {
		\Product {i \in K} {\mu_{\xi^*_i}}
	}
	\geq \Card \IClass.
\]
If $\Card \IClass = \Card {\Product {i \in K} {\mu_{\xi^*_i}}}$ we can
define \CritInd to be $K$. Otherwise, our assumption on the cardinal
arithmetic gives
\[
	\Card[\Big] {
		\Product {i \in K} {\mu_{\xi^*_i}}
	}
	= \SuccCard{\Par[\Big]{
		\BigUnion {i \in K} {\mu_{\xi^*_i}}
	}}
	> \Card \IClass.
\]
By the choice of $\alpha^*$, $\Card \IClass \geq \sup \SmallCards \geq
\BigUnion {i \in K} {\mu_{\xi^*_i}}$. Hence $\Card \IClass = \BigUnion
{i \in K} {\mu_{\xi^*_i}}$ and again we can choose \CritInd to be $K$.

It follows, that when $\Card \CritJ < \kappa$ we can find \CritInd
satisfying \EquationOnPage {small_CritInd}.
%

%
\subsection{\protect %
Case 2: The parameter depends on $\kappa$-many coordinates}
\label{large}%

\Remark If \MuSeq is such that each $\mu_\xi$ is \SuccCard \kappa or
$\lambda$, we have so far proved that $\theta$ must be either $\leq
\SuccCard \kappa$ or $\theta = \lambda$.

For the rest of the proof we assume that the set \CritJ, given in
\EquationOnPage {CritJ}, has cardinality $\kappa$ and for every $K \in
\Subsets \CritJ {<\kappa}$, $\Card [\big] { \Product {i \in K}
{\mu_{\xi^*_i}} } < \Card \IClass$. So $\chi^* \leq \Card \IClass$
holds, where
\[
	\chi^* =
	\BigUnion {K \in \Subsets \CritJ {<\kappa}}
	{\Card[\Big]{ \Product{i \in K}{\mu_{\xi^*_i}} }}.
\]
As we already know that the inequality $\Card \IClass \leq \Card
[\big] {\Product {i \in \CritJ} {\mu_{\xi^*_i}}}$ holds, the remaining
problem is that why is $\theta$ a product of strictly less than
$\kappa$-many cardinals in \Set {\mu_{\xi^*_i}} {i \in \CritJ}?

\begin{DEFINITION}{enum}%
Define \ESeqs to be the set of all sequences $\EpsilonSeq- = \Seq
{\epsilon_i} {i < \kappa}$ such that
\begin{itemize}

\item for each $i \in \CritJ$, $\epsilon_i \in \MuCard i \Minus \Ps
\ParamName {\xi^*_i}$,

\item for each $i \in \kappa \Minus \CritJ$, $\epsilon_i =
\delta^*_i$, and

\item for every $i < j < \kappa$, $\EpsilonPair - i \not= \EpsilonPair
- j$.

\end{itemize}

Again, to simplify our notation, we write $\pi(\DeltaSeq-)$ for the
sequence \Seq {\Map \pi {\xi^*_i} {\delta_i}} {i < \kappa} when
\DeltaSeq- is in \ESeqs and $\pi$ in \Maps \ZSeq satisfies that $\Set
{\Pair {\xi^*_i} {\delta_i}} {i < \kappa} \Subset \Dom \pi$,

Every sequence \EpsilonSeq- in \ESeqs determines a \Forcing \ZSeq-name
\EName {\EpsilonSeq-} for a function from $\kappa$ into 2. Namely, we
define \EName {\EpsilonSeq-} to be the name $\pi(\sigma_{\alpha^*})$
where $\pi$ is the mapping in \Maps \ZSeq satisfying that $\Dom \pi =
\Set {\Pair {\xi^*_i} {\delta^*_i}} {i < \kappa}$ and $\pi
(\DeltaSeq*) = \EpsilonSeq-$.

A pair \Pair {\DeltaSeq-} {\EpsilonSeq-} of sequences in \ESeqs is
called a neat pair if for all $i < j < \kappa$, $\DeltaPair - i \not=
\EpsilonPair - j$.

Denote the set \Set {i \in \CritJ} {\delta_i = \epsilon_i}, for
$\DeltaSeq-, \EpsilonSeq- \in \ESeqs$, by \Agr {\DeltaSeq-}
{\EpsilonSeq-}.
\end{DEFINITION}

The sequence \DeltaSeq \alpha is in \ESeqs when $\alpha \in
\IClass$. Also \EName {\DeltaSeq \alpha} is the name $\sigma_\alpha$
for every $\alpha \in \IClass$.  In fact, \Set {\EName {\EpsilonSeq-}}
{\EpsilonSeq- \in \ESeqs} is the collection of all the \Forcing
\ZSeq-names which are ``isomorphic'' to the fixed representative
$\sigma_{\alpha^*}$. The reason why we introduced ``neat pairs of
sequences in \ESeqs'' is that those names, determined by sequences in
a neat pair, can be ``coherently moved'' around by automorphisms of
\Forcing \ZSeq as follows.

\begin{LEMMA}{neat_pair}%
Suppose $\DeltaSeq1, \DeltaSeq2, \EpsilonSeq1, \EpsilonSeq2 \in
\ESeqs$ are such that both \Pair {\DeltaSeq1} {\EpsilonSeq1} and \Pair
{\DeltaSeq2} {\EpsilonSeq2} are neat, and moreover, $\Agr {\DeltaSeq1}
{\EpsilonSeq1} = \Agr {\DeltaSeq2} {\EpsilonSeq2}$ holds. Then there
is an automorphism \Autom \pi of \Forcing \ZSeq such that $\Autom \pi
(\ParamName) = \ParamName$, $\Autom \pi(\EName{\DeltaSeq1}) =
\EName{\DeltaSeq2}$ and $\Autom \pi(\EName{\EpsilonSeq1}) =
\EName{\EpsilonSeq2}$. Hence for every $p \in \Forcing \ZSeq$,
\[
	p \Forces \ZSeq
		\EquivENames {\DeltaSeq1} {\EpsilonSeq1}
	\Iff
	\Autom \pi(p) \Forces \ZSeq
		\EquivENames {\DeltaSeq2} {\EpsilonSeq2}.
\]

\begin{Proof}%
There is a mapping $\pi$ in \Maps \ZSeq such that $\pi(\DeltaSeq1) =
\DeltaSeq2$ and $\pi(\EpsilonSeq1) = \EpsilonSeq2$, because the
sequences in \ESeqs are without repetition, both of the pairs are
neat, and the equation $\Agr {\DeltaSeq1} {\EpsilonSeq1} = \Agr
{\DeltaSeq2} {\EpsilonSeq2}$ holds. Furthermore, $\pi$ can be chosen
so that \Res \pi {\Ps \ParamName -} is identity and each \Map \pi
{\xi^*_i} - is a permutation of $z_{\xi^*_i}$. Hence $\pi$ determines
an automorphism as wanted.
\end{Proof}%
\end{LEMMA}

For technical reasons we define
\ARRAY[ll]{
	\SigmaAgr =
	\Set [\Big] {I \Subset \kappa}
	{& \ThereAre \alpha \not= \beta \in \IClass \SuchThat\\
	& \Pair {\DeltaSeq \alpha} {\DeltaSeq \beta}
	\Text{is neat and} 
	I \Subset \Agr {\DeltaSeq \alpha} {\DeltaSeq \beta}}.
}
The next lemma explains why we closed the set \SigmaAgr under subsets:
all the names $\sigma_\alpha$, $\alpha \in \IClass$, are forced to be
nonequivalent, and moreover, all those names are forced to be
nonequivalent, which are determined by a neat pair of sequences
agreeing in a smaller set than some pair of the fixed sequences
\DeltaSeq \alpha, $\alpha \in \IClass$.

\begin{LEMMA}{One_forces_nonequivalent}%
For all $\DeltaSeq-, \EpsilonSeq- \in \ESeqs$, if \Pair {\DeltaSeq-}
{\EpsilonSeq-} is neat and \Agr {\DeltaSeq-} {\EpsilonSeq-} is in
\SigmaAgr, then $\One \Forces \ZSeq \NotEquivENames {\DeltaSeq-}
{\EpsilonSeq-}$.

\begin{Proof}%
First we fix $\alpha \not= \beta \in \IClass$ and $I$ such that \Pair
{\DeltaSeq \alpha} {\DeltaSeq \beta} is neat and $I = \Agr
{\DeltaSeq-} {\EpsilonSeq-} \Subset \Agr {\DeltaSeq \alpha} {\DeltaSeq
\beta}$. Let \DeltaSeq' be a sequence in \ESeqs which satisfies that
\EqualRes {\DeltaSeq'} {\DeltaSeq \alpha} I and for all $i \in \CritJ
\Minus I$, $\delta'_i \not\in \Set {\delta^\alpha_j} {j <
\kappa}$. Then the pair \Pair {\DeltaSeq'} {\DeltaSeq \alpha} is neat
and $\Agr {\DeltaSeq'} {\DeltaSeq \alpha} = I$. We want to show that
$\One \Forces \ZSeq \NotEquivENames {\DeltaSeq'} {\DeltaSeq \alpha}$,
because then it follows from \Lemma {neat_pair} that $\One \Forces
\ZSeq \NotEquivENames {\DeltaSeq-} {\EpsilonSeq-}$.

Suppose, contrary to this claim, that $p \in \Forcing \ZSeq$ satisfies
\[
	p \Forces \ZSeq
	\EquivENames {\DeltaSeq'} {\DeltaSeq \alpha}.
\]
Let $J$ denote the set \Agr {\DeltaSeq \alpha} {\DeltaSeq \beta} and
choose a sequence \EpsilonSeq' from \ESeqs so that \EqualRes
{\DeltaSeq \beta} {\EpsilonSeq'} J and for all $i \in \CritJ \Minus
J$,
\[
	\epsilon'_i \not\in
	\Ps p {\xi^*_i} \Union
	\Set {\delta'_j} {j < \kappa} \Union
	\Set {\delta^\alpha_j} {j < \kappa}.
\]
Then the pair \Pair {\EpsilonSeq'} {\DeltaSeq \alpha} is neat and
$\Agr {\EpsilonSeq'} {\DeltaSeq \alpha} = J$. By the choice of the
names $\sigma_\alpha$ and $\sigma_\beta$, i.e., by \EquationOnPage
[\EquationRefCtrl] {noneq_names}, $\One \Forces \ZSeq \sigma_\alpha
\NotEquiv \sigma_\beta$. Once more, it follows from \Lemma {neat_pair}
that
\[
	\One \Forces \ZSeq
	\sigma_\alpha = \EName {\DeltaSeq \alpha}
	\NotEquiv \EName {\EpsilonSeq'}.
\]

Choose $\pi$ from \Maps \ZSeq so that $\pi(\ParamName) = \ParamName$,
$\pi(\DeltaSeq') = \DeltaSeq'$, $\pi(\DeltaSeq \alpha) =
\EpsilonSeq'$, \Res \pi {\Par [\big] {\pi[\Ps p -] \Inter \Ps p -}} is
identity, and $\pi$ determines an automorphism \Autom \pi of \Forcing
\ZSeq. This is possible by the choice of the sequence \EpsilonSeq'.
Since $\Agr {\DeltaSeq'} {\EpsilonSeq'} = \Agr {\DeltaSeq'} {\DeltaSeq
\alpha}$ and \Pair {\DeltaSeq'} {\EpsilonSeq'} is a neat pair, it
follows from \Lemma {neat_pair} that
\[
	\Autom \pi(p) \Forces \ZSeq
	\EquivENames {\DeltaSeq'} {\EpsilonSeq'}.
\]
Now there is $q \in \Forcing \ZSeq$ satisfying $q \leq p$ and $q \leq
\Autom \pi(p)$. Since \EquivName is a name for an equivalence
relation, $q \Forces \ZSeq \EquivENames {\DeltaSeq \alpha}
{\EpsilonSeq'}$, a contradiction.
\end{Proof}
\end{LEMMA}

Next we want to show that there is always a small set of indices
outside of \SigmaAgr.

\begin{LEMMA}{EqAndSmallAgr}%
When \CritJ has cardinality $\kappa$ there are $p \in \Forcing \ZSeq$
and a neat pair \Pair {\DeltaSeq-} {\EpsilonSeq-} of sequences in
\ESeqs satisfying that%
\[%
	\Agr {\DeltaSeq-} {\EpsilonSeq-} \in \Subsets \CritJ {<\kappa}
	\And
	p \Forces \ZSeq \EquivENames {\DeltaSeq-} {\EpsilonSeq-}.
\]
%

%
\begin{Proof}%
First of all, for every $i \in \CritJ$ and $\eta \in \Functions i 2$
we fix an ordinal $\beta_\eta$ from $\mu_{\xi^*_i} \Minus \Ps
\ParamName {\xi^*_i}$ so that for all $i,j \in \CritJ$, $\eta \in
\Functions i 2$, and $\nu \in \Functions j 2$, $\beta_\eta =
\beta_\nu$ iff $i = j$ and $\eta = \nu$. Fix also a coordinate $\zeta
< \lambda$ so that $\mu_\zeta > \theta$ and $\zeta \not\in Z$ \Note
{$\zeta$ is outside of \Forcing \ZSeq}. Suppose $G$ is a
$P_{\mu_\zeta}$-generic set over $V$.  For any function $u \in
\OfGenExt {\Functions \kappa 2} G$, let \DeltaSeq u denote the
following sequence: $\DeltaSeq u = \Seq {\delta^u_i} {i < \kappa}$,
$\delta^u_i = \beta_{\Res u i}$ if $i \in \CritJ$, and $\delta^u_i =
\delta^*_i$ otherwise. Then each of the sequences \DeltaSeq u is in
$\OfGenExt \ESeqs G$. Moreover, \Pair {\DeltaSeq u} {\DeltaSeq v} is a
neat pair of for all $u$ and $v$ in \OfGenExt {\Functions \kappa 2} G.

Let $H$ be a \Forcing \ZSeq-generic set over \GenExt G. In $V[G]$,
there are at least $\mu_\zeta$ many different functions from $\kappa$
into 2. By the assumption \EquationOnPage [\EquationRefCtrl]
{initial_assumption} and \Conclusion {Absolute}, there are only
$\theta$-many equivalence classes of \Equiv in \DoubleGenExt G H. It
follows, that for some $p \in H$ and $u \not= v \in \OfGenExt
{\Functions \kappa 2} G$ the following holds in \GenExt G,
\[
	p \Forces \ZSeq
	\tau_{\DeltaSeq u} \EquivName \tau_{\DeltaSeq v}.
\]
By the definition of the ordinals $\beta_\nu$, we have that $\Agr
{\DeltaSeq u} {\DeltaSeq v} = \Set {i \in \CritJ} {\EqualRes u v i}
\in \Subsets \CritJ {<\kappa}$, and hence \Agr {\DeltaSeq u}
{\DeltaSeq v} is in $V$.

Now, in $V$, we can fix a neat pair \Pair {\EpsilonSeq 1} {\EpsilonSeq
2} of sequences in \ESeqs such that $\Agr {\EpsilonSeq 1} {\EpsilonSeq
2} = \Agr {\DeltaSeq u} {\DeltaSeq v}$. Let $\pi \in \OfGenExt {\Maps
\ZSeq} G$ be such that it determines, in \GenExt G, an automorphism
\Autom \pi of \Forcing \ZSeq satisfying $\Autom \pi (\ParamName) =
\ParamName$, $\Autom \pi (\DeltaSeq u) = \EpsilonSeq 1$, and $\Autom
\pi (\DeltaSeq v) = \EpsilonSeq 2$. For such $\pi$ in \GenExt G, we
have that $\Autom \pi (p) \Forces \ZSeq \tau_{\EpsilonSeq 1}
\EquivName \tau_{\EpsilonSeq 2}$. Note, that the condition $q = \Autom
\pi (p)$ is in $V$. From the equivalence of the forcings
$P_{\mu_\zeta} \times \Forcing \ZSeq$ and $\Forcing \ZSeq \times
P_{\mu_\zeta}$, together with \ItemOfLemma {Absolute} {Pi^1_1}, it
follows that already in $V$,
\[
	q \Forces \ZSeq
	\tau_{\EpsilonSeq 1} \EquivName \tau_{\EpsilonSeq 2}.
\]
\begin{EasyComment}%
\begin{EasyNote}%
Suppose there is $r \leq q$ with
$r \Forces \ZSeq \tau_{\EpsilonSeq 1} \NotEquivName \tau_{\EpsilonSeq
2}$. By \ItemOfLemma {Absolute} {Pi^1_1}, it would be the case that
$r \Forces \MuSeq \tau_{\EpsilonSeq 1} \NotEquivName \tau_{\EpsilonSeq
2}$, and hence, also
$r \Forces {\Concat {\SimpleSeq {\mu_\zeta}} \ZSeq} \tau_{\EpsilonSeq
1} \NotEquivName \tau_{\EpsilonSeq 2}$, a contradiction.
%
%
%
%
Suppose, contrary to the claim, that in $V$, $q' \in P_{\mu_\zeta}$
and $q' {\Forces-}_{P_{\mu_\zeta}} \Quote {q \Forces \ZSeq
\tau_{\EpsilonSeq 1} \EquivName \tau_{\EpsilonSeq 2}}$. In addition,
assume that $r \leq q$ is a condition in \Forcing \ZSeq with $r
\Forces \ZSeq \tau_{\EpsilonSeq 1} \NotEquivName \tau_{\EpsilonSeq
2}$.

Let $K$ be a \Forcing \ZSeq-generic set over $V$ with $r \in K$. Then
in \GenExt K, $f_1 \NotEquiv f_2$ where $f_1$ and $f_2$ are
interpretations of the names $\tau_{\EpsilonSeq 1}$ and
$\tau_{\EpsilonSeq 2}$ over $K$, respectively. Let $H$ be a
$P_{\mu_\zeta}$-generic set over \GenExt K with $q' \in H$. Then in
\DoubleGenExt K H, still $f_1 \NotEquiv f_2$ by \Conclusion
{Absolute}.

However, $H$ is $P_{\mu_\zeta}$-generic set over $V$, too, $K$ is
\Forcing \ZSeq-generic over \GenExt H, too, and furthermore,
$\DoubleGenExt K H = \DoubleGenExt H K$, \cite [Theorem 1.4 on page
253] {Kunen}. Therefore, $q' \in H$ implies that $\GenExtSat H {\Quote
{q \Forces \ZSeq \tau_{\EpsilonSeq 1} \EquivName \tau_{\EpsilonSeq
2}}}$, and $r \in K$ implies that $q \in K$ and so \DoubleGenExtSat H
K {f_1 \Equiv f_2}, where $f_1$ and $f_2$ are the same interpretations
of $\tau_{\EpsilonSeq 1}$ and $\tau_{\EpsilonSeq 2}$ over $K$ as
above, a contradiction.
\end{EasyNote}
\end{EasyComment}
\end{Proof}%
%
%
\end{LEMMA}

Finally, we claim that $\Card \IClass = \chi^*$, and thus we can
satisfy \EquationOnPage {large_CritInd}. Suppose, contrary to this
claim, that $\Card \IClass > \chi^*$.  In the lemma below, we show
that then all the subsets of \CritJ of cardinality $< \kappa$ are in
\SigmaAgr. It follows from \Lemma {One_forces_nonequivalent}, that for
all $\DeltaSeq-, \EpsilonSeq- \in \ESeqs$, if \Pair {\DeltaSeq-}
{\EpsilonSeq-} is neat and \Agr {\DeltaSeq-} {\EpsilonSeq-} is of
cardinality $< \kappa$, then $\One \Forces \ZSeq \NotEquivENames
{\DeltaSeq-} {\EpsilonSeq-}$.  By \Lemma {EqAndSmallAgr}, this leads
to a contradiction. So it remains to prove the following last lemma.

\begin{LEMMA}{small_sets_of_CritJ_in_SigmaAgr}%
If $\Card \IClass > \chi^*$ then $\Subsets \CritJ {< \kappa} \Subset
\SigmaAgr$.

\begin{Proof}%
Fix a set $K$ from \Subsets \CritJ {< \kappa}. Since
\[
	\Card \IClass >
	\chi^* \geq
	\Card[\Big]{ \Product{i \in K}{\mu_{\xi^*_i}} } \geq
	2^\kappa,
\]
there is $X_1 \Subset \IClass$ of cardinality $\SuccCard {(2^\kappa)}$
such that for all $\alpha \not= \beta \in X_1$, $K \Subset \Agr
{\DeltaSeq \alpha} {\DeltaSeq \beta}$.  By $\Delta$-lemma one can find
$X_2 \in \Subsets {X_1} {\SuccCard {(2^\kappa)}}$ such that for all
$\alpha \not= \beta \in X_2$, the intersection $\Set {\delta^\alpha_i}
{i < \kappa} \Inter \Set {\delta^\beta_i} {i < \kappa}$ is some fixed
set $\Xi$.  There are also $I \Subset \kappa$ and $X_3 \in \Subsets
{X_2} {\SuccCard {(2^\kappa)}}$ such that for all $\alpha \in X_3$,
$\Set {i < \kappa} {\delta^\alpha_i \in \Xi} = I$.  Hence there is
$\alpha \not= \beta \in X_3$ with \EqualRes {\DeltaSeq \alpha}
{\DeltaSeq \beta} I and $\Set {\delta^\alpha_i} {i \in \kappa \Minus
I} \Inter \Set {\delta^\beta_i} {i \in \kappa \Minus I} = \emptyset$,
i.e., \Pair {\DeltaSeq \alpha} {\DeltaSeq \beta} forms a neat pair
with $K \Subset I = \Agr {\DeltaSeq \alpha} {\DeltaSeq \beta}$.
\end{Proof}
\end{LEMMA}%
%
%
\end{SECTION}

\begin{SECTION} {-} {Remarks}{Remarks}
%
%
In fact, we needed the assumption that \Equiv is
$\Sigma^1_1$-definable over \Her \kappa only in the proof of the
absoluteness of \Equiv for extensions over the subforcing \Forcing
\ZSeq and the whole forcing \Forcing \MuSeq, i.e., in the proof of
\Lemma {Absolute}. From \Conclusion {Absolute} it follows that
\Theorem {No} holds for all equivalence relations which are definable
over \Her \kappa using a parameter and a sentence which is a Boolean
combination of a sentence containing one second order existential
quantifier \Note {$\Sigma^1_1$-sentence} and a sentence containing one
second order universal quantifier \Note {$\Pi^1_1$-sentence}. This
observation has a minor application in \cite {ShVaWeakly}.  Note, that
there is in preparation by Shelah a continuation of this paper where
the result is generalized \Note {for example the singular cardinal
hypothesis will be eliminated}. For a more general treatment of the
subject see \cite {Sh664}.

Note that the possible numbers of $\kappa$-branches in trees of
cardinality $\kappa$ and the possible numbers of equivalence classes
of $\Sigma^1_1$-equivalence relations are consistent wise almost the
same. The main difference is of course the number \SuccCard \kappa
\Note {and 0, too}. Particularly, if $\ChiSeq = \Seq {\chi_i} {i <
\gamma}$ is a sequence of nonzero cardinals such that $\gamma \leq
\kappa$ and for every $i < \gamma$, there exists a tree $T_i$ with
$\Card {T_i} \leq \kappa$ and $\Card {\Br \kappa {T_i}} = \chi_i$,
then
there exists a tree $T$ with $\Card T \leq \kappa$ and $\Card {\Br
\kappa T} = \BigUnion {i < \gamma} {\chi_i}$, and furthermore,
there exists a tree $T$ with $\Card T \leq \kappa$ and $\Card {\Br
\kappa T} = \Card {\Product {i < \gamma} {\chi_i}}$, provided that
$\gamma < \kappa$ and $\kappa^{< \kappa} = \kappa$.
Therefore, a variant of the \Theorem {No} holds, where instead of the
possible numbers of equivalence classes one considers the numbers of
$\kappa$-branches through trees of cardinality $\kappa$.



The following facts are also useful to know, when applying the theorem
proved.  Write \Fn \kappa \kappa for the ordinary Cohen-forcing which
adds a generic subset of $\kappa$, i.e., the forcing
\[
	\Set {\eta}
	{\eta \Text {
		is a partial function from $\kappa$ into 2 and
	} \Card \eta < \kappa}
\]
ordered by reverse inclusion.

\begin{FACT}{locally_Cohen}%
\begin{ITEMS}%
\ITEM{one_P}%
There is a dense subset $Q \Subset \Fn \kappa \kappa$ and a dense
embedding of $Q$ into $P_\kappa$ \Note {where $P_\kappa$ is the
forcing adding a tree with $\kappa$-many branches, see \Definition
{TForc}}.

\ITEM{equivalent}%
Every subforcing \Forcing \ZSeq of \Forcing \MuSeq is equivalent to
\Fn \kappa \kappa provided that the length of \ZSeq is at most
$\kappa$ and each $z_\xi$ has cardinality $\kappa$.

\ITEM{locally_Cohen}%
The forcing \Forcing \MuSeq is locally $\kappa$ Cohen, i.e., every
subset $Q$ of \Forcing \MuSeq of size $\leq \kappa$ is included in a
complete subforcing $Q'$ of \Forcing \MuSeq so that $Q'$ is equivalent
to \Fn \kappa \kappa.

\ITEM{preserves_wc}%
Assume that $\kappa$ is a weakly compact cardinal, and $V$ is such
that $\kappa$ remains weakly compact after forcing with \Fn \kappa
\kappa. Then every locally $\kappa$ Cohen forcing preserves weakly
compactness of $\kappa$.

\end{ITEMS}%
\end{FACT}

Note in the last claim, that if $\kappa$ is a weakly compact cardinal
then, using upward Easton forcing, it is possible to have a generic
extension \GenExt H such that $\kappa$ is weakly compact in \GenExt H
and $\kappa$ remains weakly compact in all extensions \DoubleGenExt H
G, where $G$ is \Fn \kappa \kappa-generic over \GenExt H (Silver).
\PvComment {\cite [Exercise I9 on page 298] {Kunen}}
These facts are applied in \cite {ShVaWeakly}.


%
%
\end{SECTION}

%
%

%

%

\begin{tabbing}
Saharon Shelah:\= \\
	\>Institute of Mathematics\\
	\>The Hebrew University\\
	\>Jerusalem. Israel\\
\\
	\>Rutgers University\\
	\>Hill Ctr-Busch\\
	\>New Brunswick. New Jersey 08903\\
	\>\texttt{shelah@math.rutgers.edu}
\end{tabbing}

\begin{tabbing}
Pauli V\"{a}is\"{a}nen:\= \\
	\>Department of Mathematics\\
	\>P.O. Box 4\\
	\>00014 University of Helsinki\\
	\>Finland\\
	\>\texttt{pauli.vaisanen@helsinki.fi}
\end{tabbing}

%

\end{document}